\documentclass[journal]{new-aiaa} 
\usepackage[utf8]{inputenc}

\usepackage{algorithm}
\usepackage{algorithmic}

\usepackage{graphicx}
\usepackage{amsmath}
\usepackage{mathtools}
\usepackage[version=4]{mhchem}
\usepackage{siunitx}
\usepackage{longtable,tabularx}
\usepackage{bm}

\usepackage{subfig}
\usepackage{multicol}
\usepackage{multirow}
 
\newcommand{\x}{\bm{x}}

\newcommand{\As}{{\mathcal{A}}}
\newcommand{\A}{\bm{A}}

\newcommand{\algorithmicinput}{\textbf{input}}
\newcommand{\algorithmicoutput}{\textbf{output}}

\newcommand{\INPUT}{\item[\algorithmicinput]}
\newcommand{\OUTPUT}{\item[\algorithmicoutput]}

\title{High-dimensional {\color{black}Bayesian} optimization using both random and linear embeddings}

\author{R\'emy Priem\footnote{Artificial Intelligence Engineer for Aeronautical Systems, Technical Unit, remy.priem@intradef.gouv.fr}}
\affil{DGA Maîtrise de l'Information, Bruz, France}
\author{Youssef Diouane\footnote{Professor, Mathematics and Industrial Engineering Department, youssef.diouane@polymtl.ca, AIAA MDO TC Member.}}
\affil{GERAD \& Polytechnique Montr\'eal, Montreal, QC, Canada}
\author{Nathalie Bartoli\footnote{Senior researcher, Information Processing and Systems Department, nathalie.bartoli@onera.fr, AIAA MDO TC Member.}, Sylvain Dubreuil\footnote{Researcher, Information Processing and Systems Department, sylvain.dubreuil@onera.fr, AIAA Member.}, Paul Saves\footnote{Post-doctoral fellow, Information Processing and Systems Department, paul.saves@onera.fr, Student AIAA Member.}} 
\affil{DTIS, ONERA, Université de Toulouse, Toulouse, France}
\affil{Fédération ENAC ISAE-SUPAERO ONERA, Université de Toulouse, 31000, Toulouse, France}

\begin{document}

\maketitle

\begin{abstract}
Bayesian optimization (BO) is one of the most powerful strategies to solve {\color{black}computationally expensive-to-evaluate} blackbox optimization problems.
However, BO methods are conventionally used for optimization problems of small dimension because of the curse of dimensionality.
In this paper, to solve high dimensional optimization problems, we propose to incorporate linear embedding subspaces of small dimension to efficiently perform the optimization. An adaptive learning strategy for these linear embeddings is carried out in conjunction with the optimization.
The resulting BO method, named EGORSE, combines in an adaptive way both random and supervised linear embeddings.
EGORSE has been compared to state-of-the-art algorithms and tested on academic examples with a number of design variables ranging from 10 to 600.
The obtained results show the high potential of EGORSE to solve high-dimensional blackbox optimization problems, both in terms of CPU time and {\color{black}{limited}} number of calls to the expensive blackbox {\color{black}{simulation}}.
\end{abstract}

\section{Nomenclature}

{\renewcommand\arraystretch{1.0}
\noindent\begin{longtable*}{@{}l @{\quad=\quad} l@{}}
$\Omega$ & design space \\
$d$ & dimension of the design space \\
$f$ & objective function \\
$\bm{x}$ & vector of design variables\\
$\mu$ & prior mean function \\
$k$ & covariance kernel\\
$\mathcal{D}^{(l)}$ & design of experiments of $l$ sampled points \\
$y$ & output of the objective function \\
$\mathcal{N}$ & Gaussian distribution \\
$\hat\sigma^{(l)}$ & variance function of a Gaussian process conditioned by $l$ sampled points \\
$\hat\mu^{(l)}$ & mean function of a Gaussian process conditioned by $l$ sampled points \\
$\bm{\theta}^{(l)}$ & hyper-parameters of a Gaussian process conditioned by $l$ sampled points \\
$\alpha^{(l)}$ & acquisition function \\
$y_{min}^{(l)}$ & minimum output of the objective function in a set of $l$ sampled points\\
$f_{\mathcal{A}}$ & objective function reduced in the linear subspace $\mathcal{A}$ \\
$\bm{A}$ & transfer matrix \\
$\mathcal{A}$ & linear subspace \\
$\mathcal{B}$ & hyper-cube in the linear subspace \\
$T$ & number of dimension reduction methods \\
$\mathcal{R}$ & dimension reduction method \\
$\gamma$ & reverse application \\
$\bm{u}$ & vector of design variables in the linear subspace \\
$g$ & constraint in the linear subspace \\
\end{longtable*}}

\section{Introduction}

{\color{black} The \textit{Multidisciplinary Design Analysis and Optimization} (MDAO) methodology~\cite{cramerProblemFormulationMultidisciplinary1994} is becoming a crucial pillar in many industries, especially in aircraft design. It considers a set of interacting disciplines required to design an aircraft at the beginning of the design process. For instance, one can tackle the aerodynamics, propulsion, and structural coupling by considering the design variables of these three disciplines. Motivated by the desire to improve the environmental footprint of aviation, new cutting-edge aircraft architectures are being studied to drastically decrease aircraft consumption~\cite{schmollgruber2019multidisciplinary, priemEfficientApplicationBayesian2020a, bartoliAdaptiveModelingStrategy2019, BouhlelEfficientglobaloptimization2018,Mixed_Paul_PLS}.
In this context, it is not desirable to perform MDAO procedures with low-fidelity models developed for classical tube-and-wing aircraft configurations. Aircraft designers often use high-fidelity models, even at the beginning of the process, to accurately assess the performance of the aircraft under consideration. The resulting MDAO procedures are  thus performed on computationally time-consuming tasks to find the best architecture. Furthermore, the disciplines are often blackboxes because of their complexity, meaning that no information (i.e., derivatives) other than the objective function value is available. For this reason, MDAO cannot be completed with classical gradient-based algorithms, for which the finite-difference method would be not adapted (due to the presence of noise) or too costly.  The MDAO problems for aircraft design are computationally expensive, noisy constrained, blackbox optimization problems with a large number of design variables~\cite{bartoliAdaptiveModelingStrategy2019, feliot2016design}. In this setting, evolutionary-based optimization algorithms, for which the number of function evaluations required for optimization is often too large, are not adapted as well~\cite{schmollgruber2019multidisciplinary}. }

In the last two decades, this issue is tackled with \textit{Bayesian Optimization} (BO) framework~\cite{frazierTutorialBayesianOptimization2018,JonesEfficientglobaloptimization1998,MockusBayesianmethodsseeking1975,ShahriariTakingHumanOut2016} relying on iterative enrichments of surrogate models (i.e. a \textit{Gaussian Process} (GP)) to seek for the optimum in a given design space $\Omega \subset \mathbb{R}^d$ where $d \in \mathbb{N}^+$ is the number of design variables. 
{\color{black} BO is particularly well-suited for MDAO due to its reliance on GP for the automatic enrichment of the surrogate model. 
By leveraging the predictions and uncertainties provided by GP, the enrichment process effectively targets promising areas and highly uncertain regions where the optimum may be concealed.}
Because of the characteristics of the aircraft design, like safety constraints and the important number of design variables, the resulting MDAO procedures need specific BO algorithms to be solved.
\textit{Constrained Bayesian Optimization} (CBO)~\cite{frazierTutorialBayesianOptimization2018,gelbartBayesianOptimizationUnknown2014a,ShahriariTakingHumanOut2016,bartoliAdaptiveModelingStrategy2019,YDIOUANE_2023} deals with the constraints whereas \textit{High Dimensional Bayesian Optimization} (HDBO) undertakes the high number of design variables~\cite{erikssonScalableGlobalOptimization2019a,WangBayesianoptimizationbillion2016,WangBatchedHighdimensionalBayesian2018,binoisChoiceLowdimensionalDomain2020,KandasamyHighdimensionalBayesian2015} in these problems. HDBO problems can be mathematically described as follows:
\begin{equation}
    \min\limits_{\x \in \Omega} f(\x),
    \label{eq:opt_prob}
\end{equation}
where $f : \mathbb{R}^d \mapsto \mathbb{R}$ is the objective function, the design space $\Omega = [-1,1]^d$ is a bounded domain and $d \gg 20$ corresponds to  a high number of design variables.

In the context of BO, two approaches are mainly investigated to handle the high number of design variables.
The first one is typically based on a specific adaptation of the GP structure~\cite{BouhlelEfficientglobaloptimization2018,KandasamyHighdimensionalBayesian2015,erikssonScalableConstrainedBayesian2020} to deal with the large number of design variables, e.g., EGO-KPLS~\cite{BouhlelEfficientglobaloptimization2018} where a \textit{partial least squares} (PLS) method is used to reduce the number of the GP hyper-parameters.
The second approach focuses on the use of dimension reduction methods to scale down the design space~\cite{WangBayesianoptimizationbillion2016,WangBatchedHighdimensionalBayesian2018,binoisChoiceLowdimensionalDomain2020}, e.g., REMBO~\cite{WangBayesianoptimizationbillion2016} where a random linear embedding of the initial space is used.
The dimension of the obtained subspace can be hence much lower than the original one.
Most of the methods based on the REMBO paradigm have difficulties particularly to derive the new bounds for the reduced-dimension optimization problem~\cite{WangBatchedHighdimensionalBayesian2018}.
An important computational effort in~\cite{binoisChoiceLowdimensionalDomain2020} is in general needed to complete the optimization.

Two drawbacks of most existing HDBO methods are the computational effort needed to perform the optimization process and the construction of accurate bounds (over the subspace) in which the optimization is performed. 
In this paper, a new HDBO method, named \textit{Efficient Global Optimization coupled with Random and Supervised Embedding} (EGORSE) is introduced to overcome the challenges previously detailed.
First, the standard BO framework is described in Section \ref{sec:BO}.
Then, in Section \ref{sec:EGFORSE}, the proposed methodology is detailed. A sensitivity analysis study with respect to EGORSE hyper-parameters is detailed in Section \ref{sec:sensitivity-anal}. Using academic benchmark problems, a comparison with state-of-art HDBO methods shows the high potential of EGORSE both in terms the efficiency and the global computational effort (see Sections \ref{sec:num} and \ref{sec:rover}). Conclusions and perspectives are drawn in Section \ref{sec:conclusion}.

\section{Bayesian optimization}
\label{sec:BO}

\subsection{The general framework}

To solve the unconstrained optimization problem~\eqref{eq:opt_prob}, the \textit{Bayesian optimization} (BO) framework~\cite{MockusBayesianmethodsseeking1975,JonesEfficientglobaloptimization1998} builds a surrogate model (using a Gaussian process~\cite{RasmussenGaussianprocessesmachine2006, Krigestatisticalapproachbasic1951}) of the objective function $f$ using a set of $l$ sampled points in the design domain $\Omega$, known as the \textit{Design of Experiments} (DoE).
The optimal solution is estimated by iteratively enriching a \textit{Gaussian process} (GP)~\cite{RasmussenGaussianprocessesmachine2006,Krigestatisticalapproachbasic1951} via a search strategy that balances the exploration of the design space $\Omega$ and the minimization of the surrogate model of $f$.
Namely, at each iteration of a given BO method, the search strategy requires solving a maximization sub-problem where the objective is referred as the acquisition function~\cite{frazierTutorialBayesianOptimization2018, ShahriariTakingHumanOut2016, WangMaxvalueentropysearch2017,bartoliAdaptiveModelingStrategy2019,tranAphBO2GP3BBudgetedAsynchronouslyparallel2020}.
The acquisition function being fully defined using the GP, the search strategy is computationally inexpensive and straightforward compared to the original optimization problem~\eqref{eq:opt_prob} in which the function $f$ is assumed to be expensive-to-evaluate.
The DoE is updated iteratively using the optimal solutions of the sub-problems.
The same process is repeated until a maximum number of evaluations is reached. 
The main steps of the BO framework, when applied to solve the optimization problem~\eqref{eq:opt_prob}, are summarized in Algorithm~\ref{alg:BO}.
\begin{algorithm}[ht!]
     \begin{algorithmic}[1]
        \INPUT{: Objective function, initial DoE, a maximum number of iterations max\_nb\_it\;}
        \FOR{$l = 0$ \TO \mbox{max\_nb\_it} - 1}
            \STATE {Build the surrogate model using a GP\;}
            \STATE {Find $\bm{x}^{(l+1)}$ a solution of the enrichment maximization sub-problem\;}
            \STATE {Evaluate the objective function at $\bm{x}^{(l+1)}$\;}
            \STATE {Update the DoE\;}
        \ENDFOR
        \OUTPUT{: The best point found in the DoE\;}
    \end{algorithmic}
    \caption{The Bayesian optimization framework.}
    \label{alg:BO}
\end{algorithm}

In the next two sections, we give more details on the GP and the enrichment process.

\subsection{Gaussian process}
\label{ssec:GP}
In the context of a BO process, a scalar GP~\cite{RasmussenGaussianprocessesmachine2006, Krigestatisticalapproachbasic1951} is a surrogate model whose distribution is fully described by a prior mean function, a covariance kernel and a set of sampled points.
The global behavior of the GP is depicted by the prior mean function whereas the covariance kernel characterizes the similarities between two sampled points of the design space $\Omega$.
Let a non-conditioned scalar GP defined by a prior mean function $\mu: \mathbb{R}^d \mapsto \mathbb{R}$ and a covariance kernel $k: \mathbb{R}^{d \times 2} \mapsto \mathbb{R}$.
A conditioned scalar GP of $f$ by the DoE of $l$ sampled points $\mathcal{D}^{(l)} = \left\{ \x^{(k)}, y^{(k)} \right\}_{k=0,\ldots,l-1}$, where  $\x^{(k)} \in \Omega$ and $y^{(k)} = f\left(\x^{(k)}\right)$, defines a Gaussian distribution  $\mathcal{N} \left( \hat\mu^{(l)}, \left[\hat\sigma^{(l)}\right]^2 \right)$ for each $\x \in \Omega$.
The mean $\hat\mu^{(l)} : \mathbb{R}^d \mapsto \mathbb{R}$ and standard deviation $\hat\sigma^{(l)} : \mathbb{R}^d \mapsto \mathbb{R}$ are expressed as follows:
\begin{gather}
    \label{eq:mmu}
    \hat\mu^{(l)}(\x) =  \mu(\x) + \bm{k}^{(l)}(\x)^\top \left[{\bm{K}^{(l)}}\right]^{-1} \left(\bm{Y}^{(l)} - \bm{\mu}^{(l)}\right), \\
    \label{eq:sigma}
    \hat\sigma^{(l)}(\x) = \left( k(\x,\x) - \bm{k}^{(l)}(\x)^\top \left[{\bm{K}^{(l)}}\right]^{-1} \bm{k}^{(l)}(\x) \right)^{\frac{1}{2}},
\end{gather}
where $\bm{\mu}^{(l)} = \left[\mu\left(\x^{(0)}\right), \ldots, \mu\left(\x^{(l-1)}\right) \right]^\top$ is the prior mean vector computed on the sampled points of $\mathcal{D}^{(l)}$,  \mbox{$\bm{k}^{(l)} = \left[k\left(\x, \x^{(0)}\right), \ldots, k\left( \x, \x^{(l-1)}\right) \right]^\top$} is the covariance vector between $\x$ and the sampled points of $\mathcal{D}^{(l)}$, \linebreak $\bm{K}^{(l)} = \left[k\left(\bm{x}^{(i)},\bm{x}^{(j)}\right)\right]_{i,j=0, \ldots, l-1}$ is the covariance matrix computed on the  $\mathcal{D}^{(l)}$, and $\bm{Y}^{(l)}= \left[y^{(0)}, \ldots, y^{(l-1)}\right]^\top$
is a vector of outputs of $f$.
Note that there is a wide range of prior mean functions and covariance kernels~\cite{duvenaudStructureDiscoveryNonparametric2013} and their selection is very case dependent.
Most of these functions depend on hyper-parameters, denoted by $\bm{\theta}^{(l)} \in \left[\mathbb{R}^+\right]^n$, that need to be estimated to explain the best the DoE of the objective function $f$.
To estimate the hyper-parameters $\bm{\theta}^{(l)}$ of the GP at each iteration, a maximum likelihood estimator is typically used~\cite{gelmanBayesianDataAnalysis2014}.

\subsection{The enrichment sub-problem}

The BO framework combines the information provided by the GP (namely, $\hat\mu^{(l)}$ and $\hat\sigma^{(l)}$) to build the enrichment strategy. The latter is guided by the following maximization sub-problem:
\begin{equation}
    \max\limits_{\x \in \Omega} \alpha^{(l)}(\x),
    \label{eq:iner_opt}
\end{equation}
where $\alpha^{(l)}: \mathbb{R}^d \mapsto \mathbb{R}$ is the acquisition function.
There are numerous acquisition functions in the literature~\cite{frazierTutorialBayesianOptimization2018, ShahriariTakingHumanOut2016, WangMaxvalueentropysearch2017, bartoliAdaptiveModelingStrategy2019}, their choice is essential for the enrichment process.
The \textit{Expected Improvement} (EI)~\cite{JonesEfficientglobaloptimization1998} acquisition function is the most used in BO.
Considering the $l^{\mbox{th}}$ iteration of the BO framework, the expression $\alpha_{EI}^{(l)}$ depends on the $\hat\mu^{(l)}$ and $\hat\sigma^{(l)}$.
For a given point $\x\in \Omega$, if  $\hat\sigma^{(l)}(\x) = 0$, then $\alpha_{EI}^{(l)}(\x)=0$.
Otherwise,
\begin{equation}
    \small
    \alpha_{EI}^{(l)}(\x)=\left(y_{min}^{(l)} - \hat{\mu}^{(l)}(\x)\right) \Phi \left( \frac{y_{min}^{(l)} - \hat{\mu}^{(l)}(\x)}{\hat\sigma^{(l)}(\x)} \right) + \hat\sigma^{(l)}(\x) \phi \left(  \frac{y_{min}^{(l)} - \hat{\mu}^{(l)}(\x)}{\hat\sigma^{(l)}(\x)} \right),
    \label{eq:EI}
\end{equation}
where the functions $\Phi$ and $\phi$ are, respectively, the cumulative distribution function and the probability density function of the standard normal distribution. The current minimum is given by $y_{min}^{(l)} = \min \bm{Y}^{(l)}$.
In this framework, it is possible to tackle problems with non linear constraints~\cite{frazierTutorialBayesianOptimization2018,priemOptimisationBayesienneSous2020,ShahriariTakingHumanOut2016} using different mechanisms with different computational costs.  
However, the classical BO process can not handle high dimensional problem because of the GP model.
Indeed, building a conventional GP in high-dimension can be problematic due to the likelihood maximization step used to estimate the hyper-parameters.
Moreover, classical GPs tend to miss some information in high dimension as the distance between points increase.
In the next section, we propose a solution to help overcome these challenges.

\section{Supervised linear embeddings for BO}
\label{sec:EGFORSE}
\subsection{Method description}
\label{ssec:outline}

In this paper, the objective function is supposed to depend only on the effective dimensions 
$d_e \ll d$ where we typically assume that it exists a function $f_{\As}: \mathbb{R}^{d_e} \mapsto \mathbb{R}$ such as $f_{\As}(\A\x) = f(\x)$ with $\A \in \mathbb{R}^{d_e \times d}$, $\As = \left\{ \bm{u} = \A\x \ ; \ \forall \x \in \Omega \right\}$ and\ $\Omega = [-1,1]^d$~\cite{WangBatchedHighdimensionalBayesian2018,binoisChoiceLowdimensionalDomain2020}.
The idea is then to perform the optimization procedure in the reduced linear subspace $\As$ so that the number of hyper-parameters to estimate and dimension of the design space are reduced to $d_e$ instead of  $d$.
This allows to build the inexpensive GPs and will ease the acquisition function optimization.
Using a subspace (based on $\A$) for the optimization requires finding the effective dimension of the reduced design space $\mathcal{B} \subset \mathbb{R}^{d_e}$ as well as the backward application $\gamma: \mathcal{B} \mapsto \Omega$.

The proposed method focuses on the definition of the optimization problem when a linear subspace is used as well as on efficient construction procedure of such embedding subspaces.  
Most existing HDBO methods rely on random linear subspaces  meaning that no information is used to incorporate a priori information from the optimization problem within the embedding space which may slow down the optimization process. 
In this work, a recursive search, with $T\in \mathbb{N}$ supervised reduction dimension methods, is performed to find supervised linear subspace so that the most important search directions for the exploration of the objective function are included. 
Using an initial DoE, one can use a \textit{Partial Least Squares} (PLS) regression \cite{hellandStructurePartialLeast1988} to build such linear embedding prior to the optimization process.
Furthermore, the new search design of the optimization problem (within the linear embedding subspace) is a necessary step.
Most methods rely on a classic optimization problem formulation that may limit the process performance due to very restricted new design space.
Here, once an appropriate linear subspace is found, the optimization problem is turned into a constrained optimization problem to limit the computational cost of the algorithm. The constrained optimization problem can be solved using a classical CBO method~\cite{frazierTutorialBayesianOptimization2018,priemOptimisationBayesienneSous2020,ShahriariTakingHumanOut2016,bartoliAdaptiveModelingStrategy2019,priemUpperTrustBound2020c}.

\subsection{Definition of the reduced search spaces}

To define the optimization problem in the low dimensional spaces, transfer matrices from the initial to the low dimensional spaces and optimization domains must be defined. 
The definition of these transfer matrices relies on a set of $T\in\mathbb{N}$ dimension reduction methods.
For each of the $T$ dimension reduction methods $\mathcal{R}^{(t)}$, the transfer matrix $\A^{(t)} \in \mathbb{R}^{d_e \times d}$ is build using $\mathcal{R}^{(t)}$ where $t \in \{1,\ldots,T\}$.
In this way, we propose to use supervised dimension reduction algorithms, like the PLS~\cite{hellandStructurePartialLeast1988}, to guide the optimization process through highly varying linear subspaces of the objective function.
This procedure allows to tackle the issue of~\citet{binoisChoiceLowdimensionalDomain2020,WangBatchedHighdimensionalBayesian2018} by relying on random linear subspaces defined with random Gaussian transfer matrices.
In their works, the optimization can be performed in a subspace in which the objective function is not varying, meaning the optimum of the objective function could not be discovered.
In the BO framework, the optimization process is usually performed in an hypercube. 
However $\As^{(t)}$ is not an hypercube.
As $\Omega = [-1,1]^d$, it is possible to compute $\mathcal{B}^{(t)} \subset \mathbb{R}^{d_e}$~\cite{binoisChoiceLowdimensionalDomain2020} the smallest hypercube containing all points of $\As^{(t)}$ such as 
\begin{equation}
    \mathcal{B}^{(t)} = \left[-\sum\limits_{i=1}^d\left| A^{(t)}_{1,i} \right|, \sum\limits_{i=1}^d\left| A^{(t)}_{1,i} \right| \right] \times \cdots \times \left[-\sum\limits_{i=1}^d\left| A^{(t)}_{d_e,i} \right|, \sum\limits_{i=1}^d\left| A^{(t)}_{d_e,i} \right| \right].
    \label{eq:borne_sub}
\end{equation}
Performing the optimization process in $\mathcal{B}^{(t)}$ leads to define a backward application from $\mathcal{B}^{(t)}$ to $\Omega$ to compute the objective function on the desired point. 

\subsection{The backward application}
\label{ssec:reverse}

\begin{figure}[ht!]
    \centering
    \includegraphics[width=0.5\textwidth]{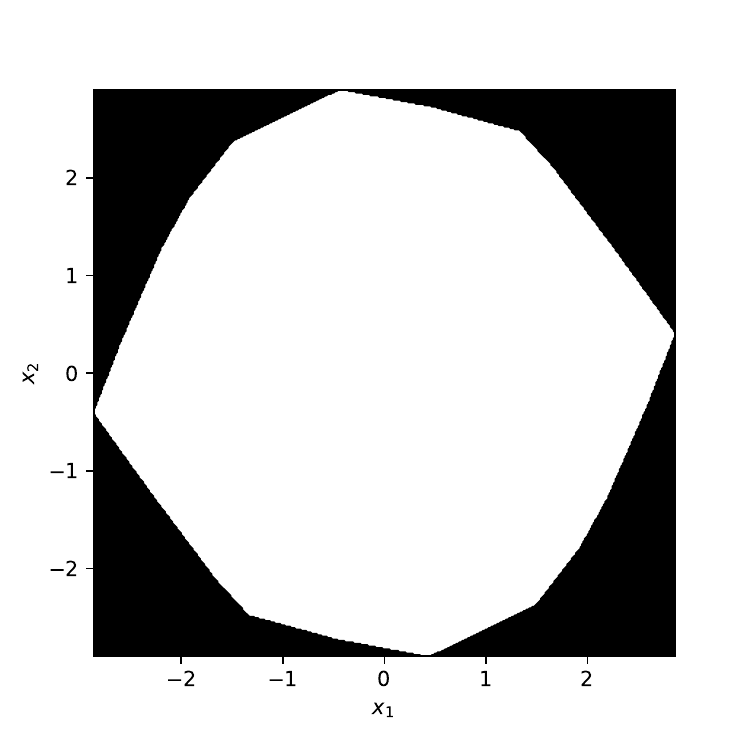}
    \caption{Illustration of $\mathcal{B}^{(t)}$ and $\As^{(t)}$. In black, the points of $\mathcal{B}^{(t)}$ without image in $\Omega$ by $\gamma_B^{(t)}$; in white, the points of $\As^{(t)}$ corresponding to the points of $\mathcal{B}^{(t)}$ with an image in $\Omega$ by $\gamma_B^{(t)}$.}
    \label{fig:backpropa}
\end{figure}

We use the backward application introduced by~\cite{binoisChoiceLowdimensionalDomain2020}. Namely, a bijective application $\gamma_B^{(t)} : \As^{(t)} \subset \mathbb{R}^d_e \mapsto \Omega \subset \mathbb{R}^d$ such that
\begin{equation}
    \gamma_B^{(t)}(\bm{u}) = \arg \min\limits_{\x \in \Omega} \left\{ \left\| \x - \left[\A^{(t)}\right]^{+} \bm{u} \right\|^2, \ \text{s.c. } \ \A^{(t)}\x = \bm{u}  \right\},
    \label{eq:quad_prob}
\end{equation}
where {$\left[\A^{(t)}\right]^{+} = \left[\A^{(t)}\right]^\top \left[ \A^{(t)} \left[\A^{(t)}\right]^\top \right]^{-1}$} the pseudo-inverse of $\A^{(t)}$. 
This problem requires to solve a quadratic optimization problem.
As $\As^{(t)} \subset \mathcal{B}^{(t)}$, $\gamma_B^{(t)}$ defines an injection from $\mathcal{B}^{(t)}$ to $\Omega$, meaning some points of $\mathcal{B}^{(t)}$ do not have any image in $\Omega$.
For instance, Figure~\ref{fig:backpropa} displays $\Omega \subset \mathbb{R}^{10}$ projected in a domain $\mathcal{B}^{(t)} \subset \mathbb{R}^2$.
The $\As^{(t)}$ domain is 
in white while the points in the black domain do not have any image in $\Omega$ by $\gamma_B^{(t)}$.

\begin{figure}[!hbt]
    \centering
    \subfloat[Projected objective function. \label{fig:proj_pb:proj}]{\includegraphics[width=0.5\textwidth]{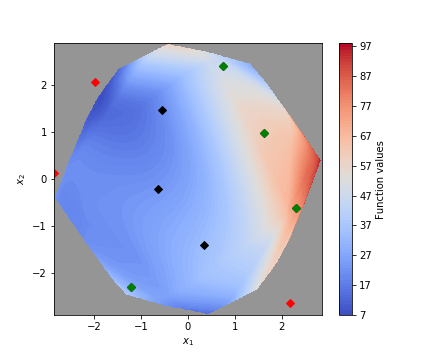}}
    \subfloat[${\alpha^{(l)}_{EI_{ext}}}$. \label{fig:proj_pb:ei}]{\includegraphics[width=0.5\textwidth]{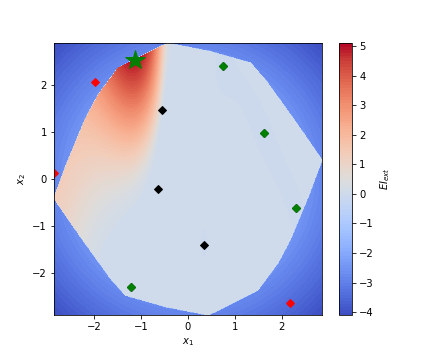}}
    \caption{Projection of an objective function of 10 design variables into a linear subspace of 2 dimensions and the corresponding ${\alpha^{(l)}_{EI_{ext}}}$ acquisition function. The grey area is the unfeasible domain, the green squares are DoE points of $\As^{(t)}$, the red squares are DoE points of $\mathcal{B}^{(t)} \backslash \As^{(t)}$ and the green star is a solution of the optimization sub-problem.}
    \label{fig:proj_pb}
\end{figure}

The method developed by~\cite{binoisChoiceLowdimensionalDomain2020} is using the $\gamma_B^{(t)}$ even if such backward application is not defined all over $\mathcal{B}^{(t)}$ by itself.
In fact, the function $f^{(t)}(\bm{u})=f\left(\gamma_B^{(t)}(\bm{u})\right)$ is only defined on $\As^{(t)}$, meaning the optimization problem in the linear subspace is given by:
\begin{equation}
    \min_{\bm{u}\in\As^{(t)}}\left\{f^{(t)}(\bm{u})\right\}.
    \label{eq:rrembo_pb}
\end{equation}
However, standard BO algorithms are only solving optimization problems with an objective function defined on an hypercube. \citet{binoisChoiceLowdimensionalDomain2020} used  $\mathcal{B}^{(t)}$, see~\eqref{eq:borne_sub}, the smallest hypercube containing $\As^{(t)}$.
This way, some points reachable by the BO algorithm might not be evaluated on the objective function as they do not have any image in $\Omega$ using $\gamma_B^{(t)}$.
To fix this issue, an additional modification of the EI acquisition function is thus introduced in \cite{binoisChoiceLowdimensionalDomain2020}. In fact, over $\As^{(t)}$, the acquisition function $\alpha^{(l)}_{EI_{ext}}(\bm{u}) = \alpha^{(l)}_{EI}(\bm{u})$~\eqref{eq:EI} while on $\mathcal{B}^{(t)} \backslash \As^{(t)}$, the acquisition function $\alpha^{(l)}_{EI_{ext}}(\bm{u}) = -\| \bm{u} \|_2$.
As the EI acquisition function is positive (i.e., lower bounded), the optimization of the acquisition function provides a point of $\As^{(t)}$. This quadratic problem~\eqref{eq:quad_prob} must be solved at each acquisition function evaluation to know if $\bm{u}$ belongs to $\As^{(t)}$.
Depending on the number of evaluations of the $\alpha^{(l)}_{EI_{ext}}$ acquisition function, an optimization process performed with the RREMBO algorithm can be expensive in CPU time.
Actually the cost in CPU time grows with the number of design variables $d$.
Figure~\ref{fig:proj_pb:proj} shows the objective function of 10 design variables projected in a linear subspace of 2 dimension.
The grey area shows that an important part of $\mathcal{B}^{(t)}$ does not have any image in $\Omega$. 
In Figure~\ref{fig:proj_pb:ei}, one can see the grey domain is replaced by negative values increasing towards the center of $\mathcal{B}^{(t)}$ allowing convergence to points having an image by $\gamma_B^{(t)}$ in $\Omega$.

\subsection{The optimization problem}
\label{ssec:sub_prob}

Section~\ref{ssec:reverse} shows that the backward application $\gamma_B^{(t)}$ is bijective from $\As^{(t)}$ to $\Omega$ and injective from $\mathcal{B}^{(t)}$ to $\Omega$. 
Thus some points of $\mathcal{B}^{(t)}$ do not have image by $\gamma_B^{(t)}$ in $\Omega$. 
To avoid points of $\mathcal{B}^{(t)} \backslash \As^{(t)}$, we define an optimization problem with a constraint which is feasible on $\As^{(t)}$ and unfeasible on $\mathcal{B}^{(t)} \backslash \As^{(t)}$.
For the value of the constraint when $\bm{u} \not\in \As^{(t)}$, we chose the opposite of the 2-norm of $\bm{u}$ .
The constraint value is negative and tends to zero when $\bm{u}$ is getting closer to $\As^{(t)}$ in norm.
To define the constraint value in the feasible zone, ones can rely on~\eqref{eq:borne_sub} and on $\Omega = [-1,1]^d$.
The equation shows that points on the edge of $\mathcal{B}^{(t)}$, having an image in $\Omega$ by $\gamma_B^{(t)}$, are in the corners of $\Omega$. 
The corners of a domain are considered to be the points whose components are $-1$ or $1$.
However, the corners of the domain are the farthest points from the center $x_c = \bm{0} \in \mathbb{R}^d$ and have the same norm. 
All the images from points of $\mathcal{B}^{(t)}$ by $\gamma_B^{(t)}$ have a lower norm than the corners of $\Omega$. 
Hence the constraint of the optimization problem is defined as the difference between the norm of the corners $\Omega$ and the norm of $\gamma_{B}^{(t)}(\bm{u})$ when $\bm{u} \in \As^{(t)}$.
Thus, the constraint value tends to zero when $\bm{u}$ is getting closer to $\mathcal{B}^{(t)} \backslash \As^{(t)}$.
Eventually, the constraint function is normalized to provide values in $[-1,1]$.
This normalization on the bounds of the constraint function balances the importance of the feasible and unfeasible domain in the optimization process.

\begin{figure}[!htb]
    \centering
    \subfloat[Optimization problem \eqref{eq:it_pb}. \label{fig:it_pb:pb}]{\includegraphics[width=0.5\textwidth]{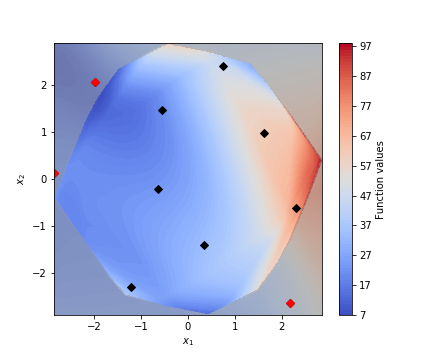}}
    \subfloat[SEGOMOE optimization sub-problem associated to problem~\eqref{eq:it_pb}. \label{fig:it_pb:spb}]{\includegraphics[width=0.5\textwidth]{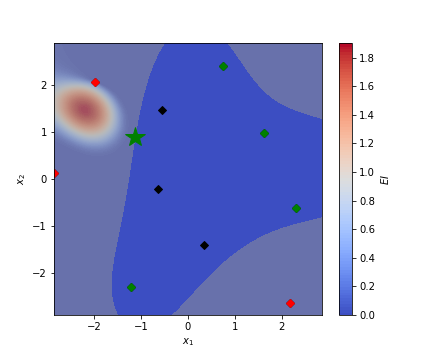}}
    \caption{An optimization problem~\eqref{eq:it_pb} in a 2 dimensional linear subspace and the associated SEGOMOE optimization sub-problem. The grey area is the unfeasible domain, the green squares are DoE points of $\As^{(t)}$, the red squares are DoE points of $\mathcal{B}^{(t)} \backslash \As^{(t)}$ and the green star is a solution of the optimization sub-problem.}
    \label{fig:it_pb}
\end{figure}

To summarize, the constraint is given by $g^{(t)}(\bm{u}) \geq 0$ where:
\begin{equation}
    g^{(t)}(\bm{u}) = \left\{
    \begin{tabular}{ll}
        $1 - \frac{\left\| \gamma_B^{(t)}\left(\bm{u}\right)\right\|^2_2}{d}$ &  if $\bm{u} \in \mathcal{A}^{(t)}$\\
        $-\left\|\bm{u}^{\A^{(t)}}\right\|^2_2$ &  otherwise
    \end{tabular}
    \right.
    \label{eq:it_cst}
\end{equation}
with $u^{\A^{(t)}}_i = u_i / \sum\limits_{j=1}^d\left| A^{(t)}_{i,j} \right|$.
The terms $u^{\A^{(t)}}_i$ and $u_i$ are respectively the $i^{\text{th}}$ components of the $\bm{u}^{\A^{(t)}}$ and $\bm{u}$ vectors.
One remarks that $f^{(t)}(\bm{u}) = f\left(\gamma_B^{(t)}(\bm{u})\right)$ function is not defined all over $\mathcal{B}^{(t)}$. 
That is why one seeks to give a value to $f^{(t)}$ on $\mathcal{B}^{(t)} \backslash \As^{(t)}$ to obtain a function not only defined on $\mathcal{B}^{(t)}$.
The extension of $f^{(t)}$ on $\mathcal{B}^{(t)} \backslash \As^{(t)}$ is given by $f^{(t)}(\bm{u}) = f\left(\gamma_W^{(t)}(\bm{u})\right)$ with  
\begin{equation}
     \gamma_W^{(t)}(\bm{u}) \in \arg \min\limits_{\x \in \Omega} \| \x - \left[\A^{(t)}\right]^{+} \bm{u} \|^2.
    \label{eq:quad_prob_wang}
\end{equation}
Indeed, the $\gamma_W^{(t)}$ application exists for all points $\bm{u} \in \mathcal{B}^{(t)}$ although it can provide the same $\x \in \Omega$ for different $\bm{u} \in \mathcal{B}^{(t)}$.
These points of $\Omega$ are moreover reachable with points of $\As^{(t)}$ using the $\gamma_B^{(t)}$ backward application.
Thus, one of the objective function minima is necessarily in $\As^{(t)}$.
Eventually, the objective function $f^{(t)}$ is given on $\mathcal{B}^{(t)}$ by 
\begin{equation}
    f^{(t)}(\bm{u}) = \left\{
    \begin{tabular}{ll}
        $f\left(\gamma_B^{(t)}\left(\bm{u}\right)\right)$    &  if $\bm{u} \in \mathcal{A}^{(t)}$\\
        $f\left(\gamma_W^{(t)}\left(\bm{u}\right)\right)$  &  otherwise
    \end{tabular}
    \right.
    \label{eq:it_obj}
\end{equation}
In fact, one solves the following constrained optimization problem in the $\mathcal{B}^{(t)}$ hypercube:
\begin{equation}
    \min\limits_{\bm{u} \in \mathcal{B}^{(t)}} \left\{ f^{(t)}(\bm{u}) \quad \text{s.c.} \quad g^{(t)}(\bm{u}) \geq 0 \right\}.
    \label{eq:it_pb}
\end{equation}
where $f^{(t)}(\bm{u})$ and $g^{(t)}(\bm{u})$ are respectively defined by Eq.~\eqref{eq:it_obj} and Eq.~\eqref{eq:it_cst}.
To solve this problem, a standard CBO algorithm~\cite{frazierTutorialBayesianOptimization2018,priemOptimisationBayesienneSous2020,ShahriariTakingHumanOut2016}, like SEGOMOE~\cite{bartoliAdaptiveModelingStrategy2019,priemUpperTrustBound2020c} can be used.
Figure~\ref{fig:it_pb:pb} shows the optimization problem~\eqref{eq:it_pb} for the $10$ design variable function given in Figure~\ref{fig:proj_pb:proj} and projected in a 2 dimensional linear subspace. 
One sees that the grey unfeasible area corresponds to the grey area of Figure~\ref{fig:proj_pb:proj}.
Figure~\ref{fig:it_pb:spb} shows the optimization sub-problem of the SEGOMOE algorithm using the EI acquisition function where the green and red squares are points of the initial DoE.
The predicted unfeasible zone, in grey in Figure~\ref{fig:it_pb:spb}, contains the unfeasible points of the initial DoE. 
Recall, with SEGOMOE, only the mean of the GP modeling the constraint is used to defined the feasible zones.
On the contrary of RREMBO no quadratic problem solving is needed to solve the optimization sub-problem. 
In fact, the introduced process avoids the quadratic problem solving in the optimization sub-problem by defining functions existing all over $\mathcal{B}^{(t)}$, which is not the case of RREMBO (see Section~\ref{ssec:reverse}).
The quadratic problem is only solved when the objective function is called at each iteration of the SEGOMOE algorithm.
Thus, the introduced EGORSE algorithm should be faster in CPU time than RREMBO.

\subsection{Adaptive learning of the linear subspace}

\begin{figure}[h]
    \centering
    \includegraphics[width=\textwidth]{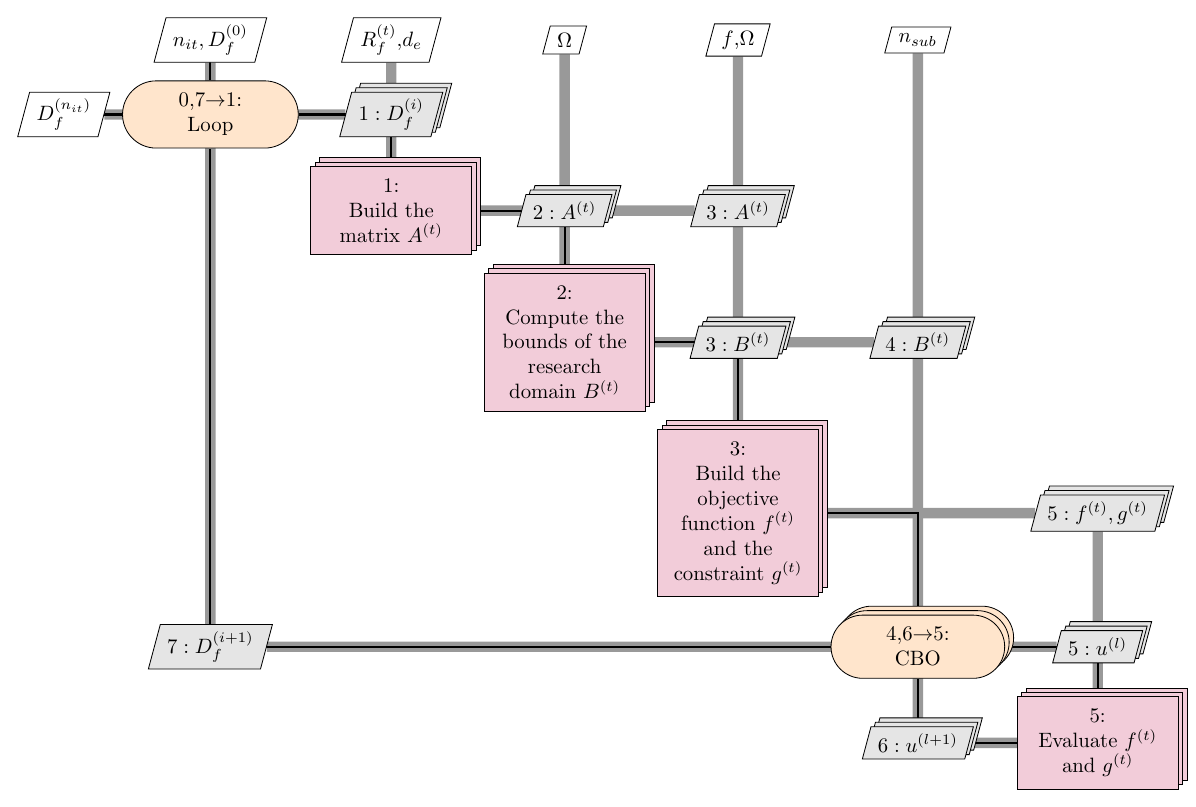}
    \caption{An XDSM diagram of the EGORSE framework.}
    \label{fig:xdsm}
\end{figure}

To ease the convergence process, we use a linear subspace discovered by supervised dimension reduction methods, i.e. taking points of the considered domain and the associated function values as inputs.
In our case, linear dimension reduction methods are used to provide the transfer matrix $\A^{(t)} \in \mathbb{R}^{d \times d_e}$ from $\Omega$ to $\As^{(t)}$.
 In order to find the minimum of the objective function, we solve the optimization problem~\eqref{eq:it_pb} with a CBO algorithm~\cite{frazierTutorialBayesianOptimization2018,ShahriariTakingHumanOut2016,priemOptimisationBayesienneSous2020} in $\mathcal{B}^{(t)}$ generated by the $\A^{(t)}$ matrix with a maximum of \texttt{max\_nb\_it\_sub} iterations.
However, in larger space, the linear subspace approximation can lack of accuracy if the number of points used by the dimension reduction methods is not sufficient.
To improve the generated subspace, the previous process is iterated using all the evaluated points during the previous iterations.
Note that points outside of the subspace $\As^{(t)}$ are added to provide information which is non biased by the subspace selection.
In that way, points coming from an optimization performed in a subspace defined by other dimension reduction methods can be used.
For instance, a transfer matrix defined by a random Gaussian distribution is chosen. 
The convergence properties, given by~\citet{binoisChoiceLowdimensionalDomain2020}, are thus preserved. 
Finally, we 
generalize this approach by using several dimension reduction methods, which can be unsupervised (like random Gaussian transfer matrices or hash tables~\cite{binoisChoiceLowdimensionalDomain2020,nayebiFrameworkBayesianOptimization2019}), or supervised (like PLS or MGP~\cite{hellandStructurePartialLeast1988,GarnettActiveLearningLinear2014}), to consider their advantages. 
The overall process of EGORSE is finally given by Algorithm~\ref{alg:EGORSE} whereas the flow chart of the method is described by the \textit{eXtended Design Structure Matrix} (XDSM) \cite{lambeExtensionsDesignStructure2012} diagram of Figure~\ref{fig:xdsm}.
\begin{algorithm}[!hbtp]
     \begin{algorithmic}[1]
        \INPUT{: an objective function $f$, an initial {DoE} $\mathcal{D}_f^{(0)}$, a maximum number of iterations \mbox{max\_nb\_it}, a maximum number of iterations by subspace max\_nb\_it\_sub, a number of active directions $d_e$, a list $\mathcal{R} = \left\{\mathcal{R}^{(1)},\ldots,\mathcal{R}^{(T)}\right\}$ of $T \in \mathbb{N}^+$ supervised or unsupervised dimension reduction methods\;}
        \FOR{$i = 0$ \TO \mbox{max\_nb\_it} - 1}
            \FOR{$t = 1$ \TO $T$}
                \STATE {Build $\A^{(t)} \in \mathbb{R}^{d_e \times d}$ with $\mathcal{R}^{(t)}$ using or not $\mathcal{D}_f^{(i)}$}
                \STATE {Define $\mathcal{B}^{(t)}$ (see Eq.~\eqref{eq:borne_sub})}
                \STATE {Build $f^{(t)}$ (see Eq.~\eqref{eq:it_obj})}
                \STATE {Build $g^{(t)}$ (seeEq.~\eqref{eq:it_cst})}
                \STATE {Solve the optimization problem
                $\min\limits_{\bm{u} \in \mathcal{B}^{(t)}} \left\{ f^{(t)}(\bm{u}) \quad \text{s.c.} \quad g^{(t)}(\bm{u}) \geq 0 \right\}$
                with a CBO algorithm using \mbox{max\_nb\_it\_sub} maximum iterations.}
            \ENDFOR
            \STATE {$\mathcal{D}_f^{(i+1)} = \mathcal{D}_f^{(i)} \cup \{\text{Points already evaluated on $f$ at iteration $i$}\}$}
        \ENDFOR
        \OUTPUT{: The best point regarding the value of $f$ in  $\mathcal{D}_f^{(\text{max\_nb\_it})}$\;}
    \end{algorithmic}
    \caption{The EGORSE process.}
    \label{alg:EGORSE}
\end{algorithm}
To conclude, the EGORSE algorithm includes three main contributions:
\begin{itemize}
    \item The restriction of the number of quadratic problems solved with a new formulation of the subspace optimization problem.
    \item The use of supervised dimension reduction methods promoting the search in highly varying directions of the objective function in $\Omega$.
    \item The adaptive learning of the linear subspace using all the evaluated points.
\end{itemize}

{\color{black}
\subsection{Efficient supervised embeddings}

Building an efficient supervised embedding (i.e., the transfer matrix $A^{(t)}$) plays a key role in the EGORSE framework. In this work, we consider two well-known supervised embedding techniques: the \textit{Partial Least-Squares (PLS)}~\cite{hellandStructurePartialLeast1988} and the \textit{Marginal Gaussian Process (MGP)}~\cite{GarnettActiveLearningLinear2014} based methods. In this section, we recall briefly the main features of the PLS and MGP methods. To simplify the notations, we will omit the iteration indices $(t)$ in this section; we will refer to the transfer matrix using within EGORSE by $A$ and to the current DoE by $\mathcal{D}_f$. 
These two embeddings PLS and MGP coupled with GP are available in the SMT toolbox~\cite{bouhlelPythonSurrogateModeling2019,SMT2023}\footnote{\href{https://github.com/SMTorg/smt}{https://github.com/SMTorg/smt}}.

\paragraph{On the PLS transfer matrix ( $\A_{\mbox{PLS}}$).}
The PLS~\cite{hellandStructurePartialLeast1988} method searches for the most important $d_e$ orthonormal directions $\bm{a}_{(i)} \in \mathbb{R}^d$, $i \in \{1,\ldots,d_e\}$ of $\Omega$ in regards of the influence of (related to the DoE $\mathcal{D}_f$) the inputs \mbox{$\bm{X}^{(l)}_{(0)} = \left[\left[\x^{(0)}\right]^\top,\ldots,\left[\x^{(l)}\right]^\top \right]^\top \in \mathbb{R}^{d \times l}$} on the outputs $\bm{Y}^{(l)}_{f,(0)} = \bm{Y}^{(l)}_f \in \mathbb{R}^l$.
These directions are recursively given by: 
\begin{gather*}
	\bm{a}'_{(i)} \in \arg \max_{\bm{a}' \in \mathbb{R}^d} \left\{ \bm{a}'^\top \left[\bm{X}^{(l)}_{(i)}\right]^\top \bm{Y}^{(l)}_{f,(i)} \left[\bm{Y}^{(l)}_{s,(i)}\right]^\top \bm{X}^{(l)}_{(i)} \bm{a}' \ , \ \bm{a}'^\top \bm{a}' = 1 \right\}, \label{eq:kpls_opt}\\
	\bm{t}_{(i)} = \bm{X}^{(l)}_{(i)} \bm{a}'_{(i)}, \quad
	\bm{p}_{(i)} = \left[\bm{X}^{(l)}_{(i)}\right]^\top \bm{t}_{(i)} , \quad c_{(i)} = \left[\bm{Y}^{(l)}_{s,(i)}\right]^\top \bm{t}_{(i)} \\ \bm{X}^{(l)}_{(i+1)} = \bm{X}^{(l)}_{(i)} - \bm{t}_{(i)} \bm{p}_{(i)}, \quad
	\bm{Y}^{(l)}_{f,(i+1)} = \bm{Y}^{(l)}_{f,(i)} - c_{(i)} \bm{t}_{(i)}.
\end{gather*}
where $\bm{X}^{(l)}_{(i+1)} \in \mathbb{R}^{d \times l}$ and $\bm{Y}^{(l)}_{f,(i+1)} \in \mathbb{R}^d$ are the residuals of the projection of $\bm{X}^{(l)}_{(i)}$ and $\bm{Y}^{(l)}_{f,(i)}$ on the $\text{i}^{\text{th}}$ principal component $\bm{t}_{(i)} \in \mathbb{R}^d$.
Finally, $\bm{p}_{(i)} \in \mathbb{R}^d$ and $c_{(i)} \in \mathbb{R}$ are respectively the regression coefficients of $\bm{X}^{(l)}_{(i)}$ and $\bm{Y}^{(l)}_{s,(i)}$ on the $\text{i}^{\text{th}}$ principal component $\bm{t}_{(i)}$ for $i \in \{1,\ldots,d_e\}$.
In fact, the square of the covariance between $\bm{a}'_{(i)}$ and $\left[\bm{X}^{(l)}_{(i)}\right]^\top \bm{Y}^{(l)}_{f,(i)}$ is recursively maximized. In this case, the PLS transfer matrix is given by:
$$
\A_{\mbox{PLS}} = \A' (\bm{P}^\top \A')^{-1},
$$
where $\A' = [\bm{a}'_{(1)},\ldots,\bm{a}'_{(d_e)}] $ and $
	\bm{P} = [\bm{p}'_{(1)},\ldots,\bm{p}'_{(d_e)}]$.

\paragraph{On the MGP transfer matrix ($\A_{\mbox{MGP}}$).}
With the MGP~\cite{GarnettActiveLearningLinear2014}, the matrix is considered as the realization of a Gaussian distribution \mbox{$\mathbb{P}(A) = \mathcal{N}\left(\A_{\mbox{MGP}}, \bm{\Sigma}_{\mbox{MGP}} \right)$} where $\A_{\mbox{MGP}}$ and $\bm{\Sigma}_{\mbox{MGP}}$ are respectively the prior mean and the covariance matrix.
The density function of $\mathbb{P}(A)$ is noted $p(\A)$.  Then, the transfer matrix $\A_{\mbox{MGP}}$ used when referring to the MGP dimension reduction method is given as the maximum of the  posterior probability law of $\A$ with respect to the DoE $\mathcal{D}_f$, i.e.,
$$
\A_{\mbox{MGP}} \in \arg \max_{A}  p\left(\A | \mathcal{D}_f\right):= p(\A) \cdot \mathcal{L}\left(\bm{Y}_f,\A\right),
$$
 where  $\mathbb{P}\left(\A | \mathcal{D}_f\right)$ represents the density function of $\mathbb{P}\left(\A | \mathcal{D}_f\right)$ and $\mathcal{L}$ the likelihood. The covariance matrix $\bm{\Sigma}_{\mbox{MGP}}$ is estimated by the inverse of the logarithm of $p\left(\A | \mathcal{D}_f\right)$ Hessian  matrix evaluated at $\A_{\mbox{MGP}}$
\begin{equation*}
            \hat{\bm{\Sigma}}^{-1}_{\mbox{MGP}}  = - \nabla^2 \log p\left(\A_{\mbox{MGP}} | \mathcal{D}_f\right),
\end{equation*}
where $\nabla^2$ is the Hessian operator with respect to $\A$. The remaining details on the MGP can be found in~\cite{GarnettActiveLearningLinear2014}.
}

\section{A sensitivity analysis over \texttt{EGORSE} hyper-parameters}
\label{sec:sensitivity-anal}
\subsection{Implementation details}

\texttt{EGORSE} method is implemented with Python 3.8. A CBO process (see Section~\ref{ssec:sub_prob}) is performed at each iteration of the \texttt{EGORSE} algorithm
using SEGO~\cite{SasenaExplorationmetamodelingsampling2002}  (from the SEGOMOE toolbox~\cite{bartoliAdaptiveModelingStrategy2019}).
The SEGOMOE toolbox uses the SMT~\cite{bouhlelPythonSurrogateModeling2019}, a Python package used to build the GP model.
In the CBO process, the EI acquisition function is optimized in two steps. 
First, a good starting point is found by solving the sub-optimization problem~\eqref{eq:iner_opt} with the ISRES~\cite{runarsson2005search} from the NLOPT~\cite{johnson2014nlopt} Python toolbox.
ISRES is an evolutionary optimization algorithm able to solve multi-modal optimization problems with equality and inequality constraints.
The algorithm explores the domain to find an optimal area maximizing the acquisition function and respecting the feasibility criteria.
However, such algorithm needs many function evaluations to converge.
To limit this number of evaluations, the solution provided by ISRES is refined with a gradient based optimization algorithm as the analytical derivatives of the acquisition function and the feasibility criteria are easily available from the GP approximations and provided as outputs from SMT~\cite{bouhlelPythonSurrogateModeling2019}.
The gradient based algorithm used is SNOPT~\cite{Gillsnopt2005} from the PyOptSparse~\cite{Perezpyopt2012} Python toolbox whose initial guess is given by ISRES.
To solve the quadratic problem~\eqref{eq:quad_prob}, the CVXOPT~\cite{andersen2013cvxopt} Python toolbox is chosen.
A point $\bm{u}\in\mathcal{B}_B^{(t)}$ is considered to belong to $\As^{(t)}$ if the CVXOPT optimization status '\verb'optimal'' is reached meaning that the quadratic problem has a solution.

\subsection{On the setting of \texttt{EGORSE} hyper-parameters}
\label{sssec:hyperparam}

\texttt{EGORSE} is controlled by a set of hyper-parameters. To select the value of these hyper-parameters, a parametric study is performed on two optimization problems.
The hyper-parameters considered in this study are the following.
\begin{itemize}
    \item \textbf{The number of points in the initial DoE}. It impacts the supervised dimension reduction methods. On the overall \texttt{EGORSE} versions, three sizes of initial DoE are tested: 5 points, $d$ points and $2d$ points where $d$ is the number of design variables.
    \item \textbf{The supervised dimension reduction method.} It changes the behavior of the algorithm favoring different directions of the domain.
    The PLS~\cite{hellandStructurePartialLeast1988} and MGP~\cite{gardnerDiscoveringExploitingAdditive2017} methods are considered in this study. 
    The PLS method is coming from the scikit-learn~\cite{scikit-learn} Python toolbox. 
    The MGP method is implemented within SMT~\cite{bouhlelPythonSurrogateModeling2019}.
    \item \textbf{The unsupervised dimension reduction method} relies on a random Gaussian transfer matrix like in~\citet{binoisChoiceLowdimensionalDomain2020} work or on hash table like in~\citet{nayebiFrameworkBayesianOptimization2019} work.
\end{itemize}
In what comes next, the following 6 possible variants of \texttt{EGORSE} will be tested and compared:
\begin{enumerate}
    \item \texttt{EGORSE Gaussian}: it uses Gaussian random transfer matrix.
    \item  \texttt{EGORSE Hash}: it uses random matrix defined by Hash table.
    \item  \texttt{EGORSE PLS}: it uses transfer matrix defined by the PLS method.
    \item  \texttt{EGORSE PLS+Gaussian}: it uses Gaussian random transfer matrix and transfer matrix defined by the PLS.
    \item  \texttt{EGORSE MGP}: it uses transfer matrix defined by the MGP method.
    \item  \texttt{EGORSE MGP+Gaussian}: it uses Gaussian random transfer matrix and transfer matrix defined by the MGP.
\end{enumerate}
\subsection{Study on two analytical problems}
\subsubsection{Definition of the two problems}

The considered class of problems is an adjustment of the Modified Branin (MB) problem~\cite{ParrInfillsamplingcriteria2012} whose number of design variables is artificially increased.
This problem is commonly used in the literature~\cite{binoisChoiceLowdimensionalDomain2020,WangBayesianoptimizationbillion2016,nayebiFrameworkBayesianOptimization2019} and it is  defined as follows:
\begin{equation}
    \min_{\bm{u} \in \Omega_1}{f_1(\bm{u})},
\end{equation}
where $\Omega_1 = [-5,10] \times [0,15]$ and
\begin{equation}
    f_1(\bm{u}) = \left[ \left( u_2 - \frac{5.1u_1^2}{4\pi^2} + \frac{5u_1}{\pi} - 6 \right)^2 +\left( 10 - \frac{10}{8\pi} \right) \cos{(u_1) + 1} \right] + \frac{5u_1 + 25}{15}.
\end{equation}
The modified version of the Branin problem is selected because it count three local minima including a global one.
The value of the global optimum is about ${\text{MB\_$d$}}_{min} = 1.1$.
Furthermore, the problem is normalized to have $\bm{u} \in [-1,1]^2$.
To artificially increase the number of design variables, a random matrix $\A_d \in \mathbb{R}^{2 \times d}$ is generated such that for all $\x \in [-1,1]^d$, $\A_d \x = \bm{u}$ belongs to $[-1,1]^2$.
An objective function MB\_$d$, where $d$ is the number of design variables, is defined such that $\text{MB\_$d$}(\x) = f_1(\A_d\x)$.
Eventually, we solve the following optimization problem:
\begin{equation}
    \min\limits_{x \in [-1,1]^d} \text{MB\_$d$}(\x) = \min\limits_{x \in [-1,1]^d} f_1(\A_d\x).
\end{equation}
In the following numerical experiments are conducted on two functions of respective dimension 10 and 100. These two test functions are denoted $\text{MB\_$10$}$ and $\text{MB\_$100$}$.

\subsubsection{Convergence plots}
\label{sssec:conv}

To study the \texttt{EGORSE} hyper-parameter impact, convergence plots are obtained for both problems.
Ten independent optimizations are thus performed on each of the problems, for each of the \texttt{EGORSE} versions, using 10 initial DoE.
The number of iterations is imposed to 10 and the number of evaluations by sub-space optimization is set to $20d_e$ and $d_e=2$. 
For  \texttt{EGORSE PLS, \texttt{EGORSE} MGP, \texttt{EGORSE} Hash and \texttt{EGORSE} Gaussian}, the number of iterations is doubled to keep a fixed number of evaluations, i.e. 800 evaluations per run.
Indeed, at each iteration, the problem is evaluated only $20d_e$ for \texttt{EGORSE} composed of a unique dimension reduction method against $2\times20d_e$ for \texttt{EGORSE} composed of two dimension reduction methods.  
All of the \texttt{EGORSE} versions are then compared displaying the evolution of the mean and standard deviation on the $10$ optimizations of the lower value of the objective function evaluated with respect to the number of evaluations. 
For readability, the standard deviation is displayed with a reduction factor of four. 

\subsubsection{Results analysis}
\label{sssec:impact}

\begin{figure}[!htbp]
    \centering
    \subfloat[DoE of 5 points.\label{fig:EGORSE_MB10_5}]{\includegraphics[width=0.50\textwidth]{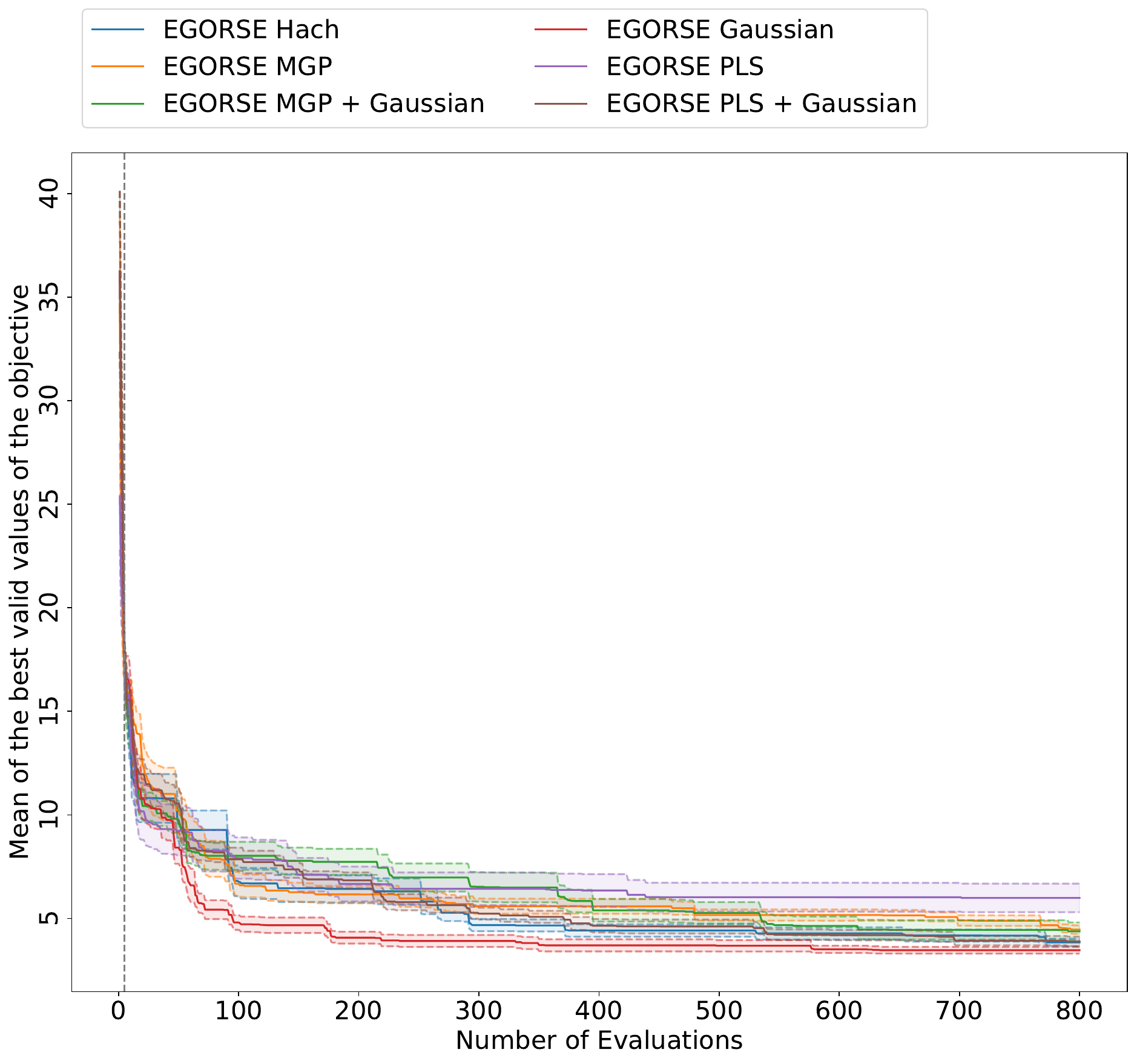}}
    \subfloat[DoE of 10 points.\label{fig:EGORSE_MB10_d}]{\includegraphics[width=0.50\textwidth]{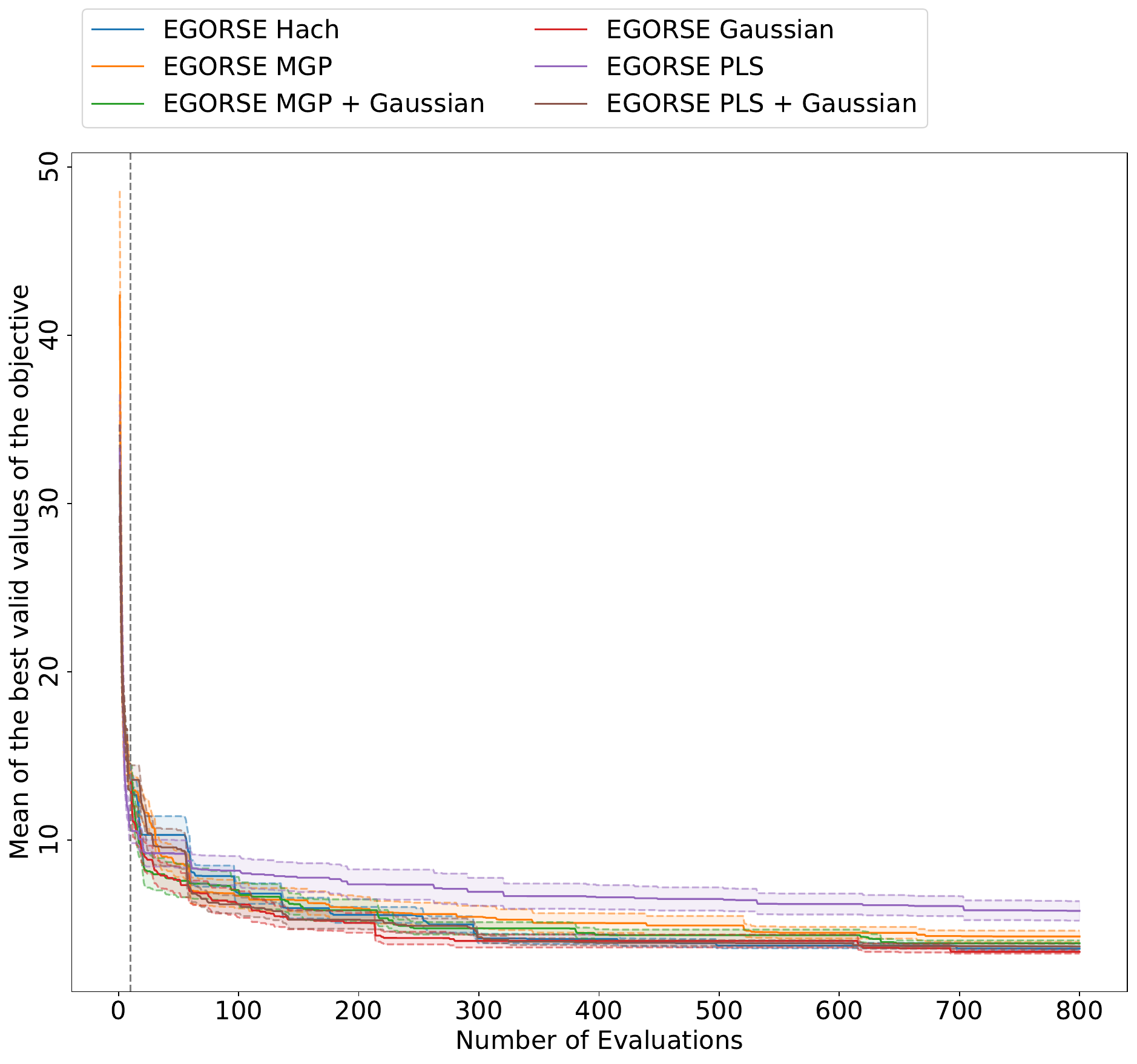}} \\
    \subfloat[DoE of 20 points.\label{fig:EGORSE_MB10_2d}]{\includegraphics[width=0.50\textwidth]{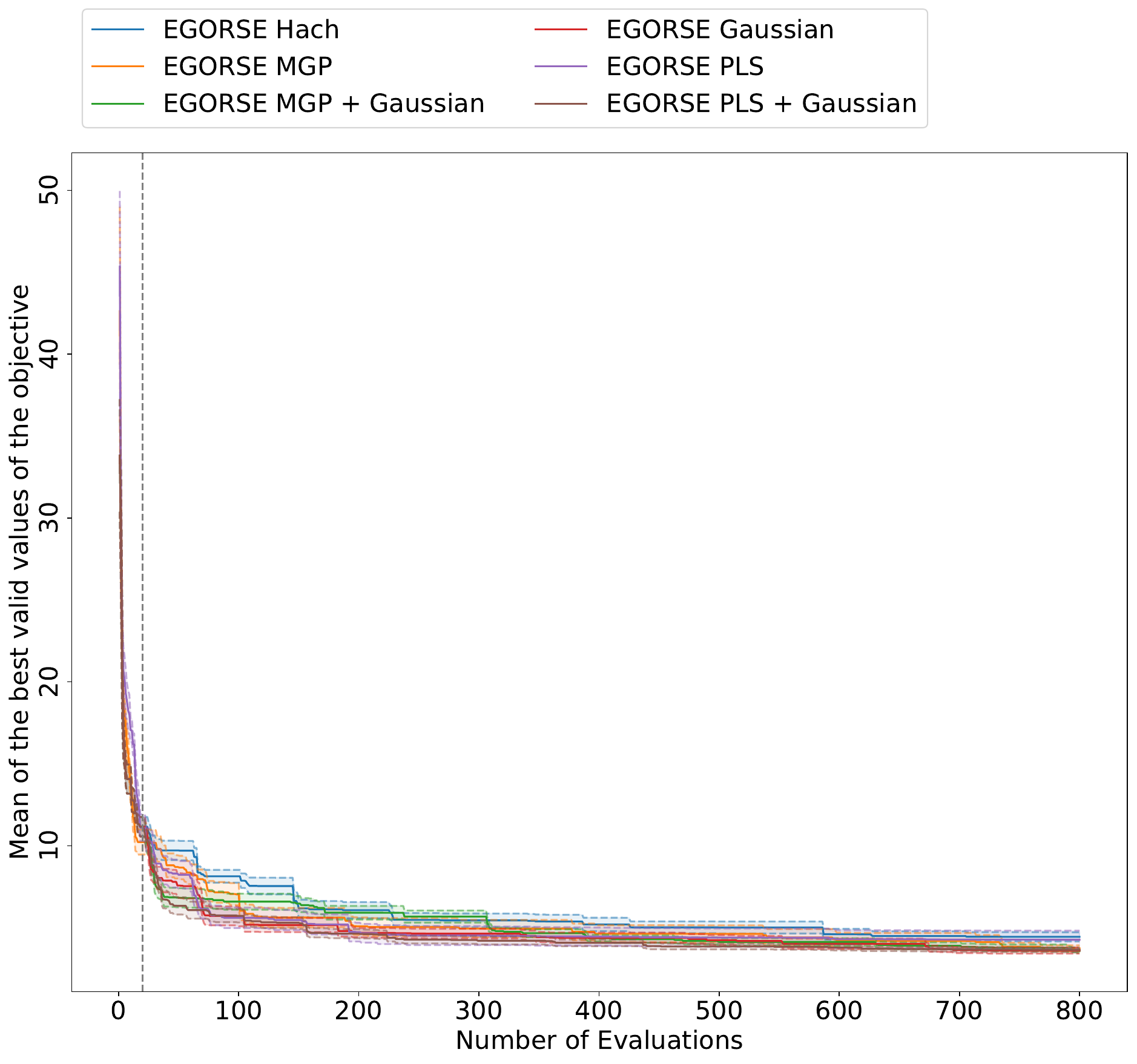}}
    \subfloat[Best versions for each DoE size. \label{fig:EGORSE_MB10_best}]{\includegraphics[width=0.50\textwidth]{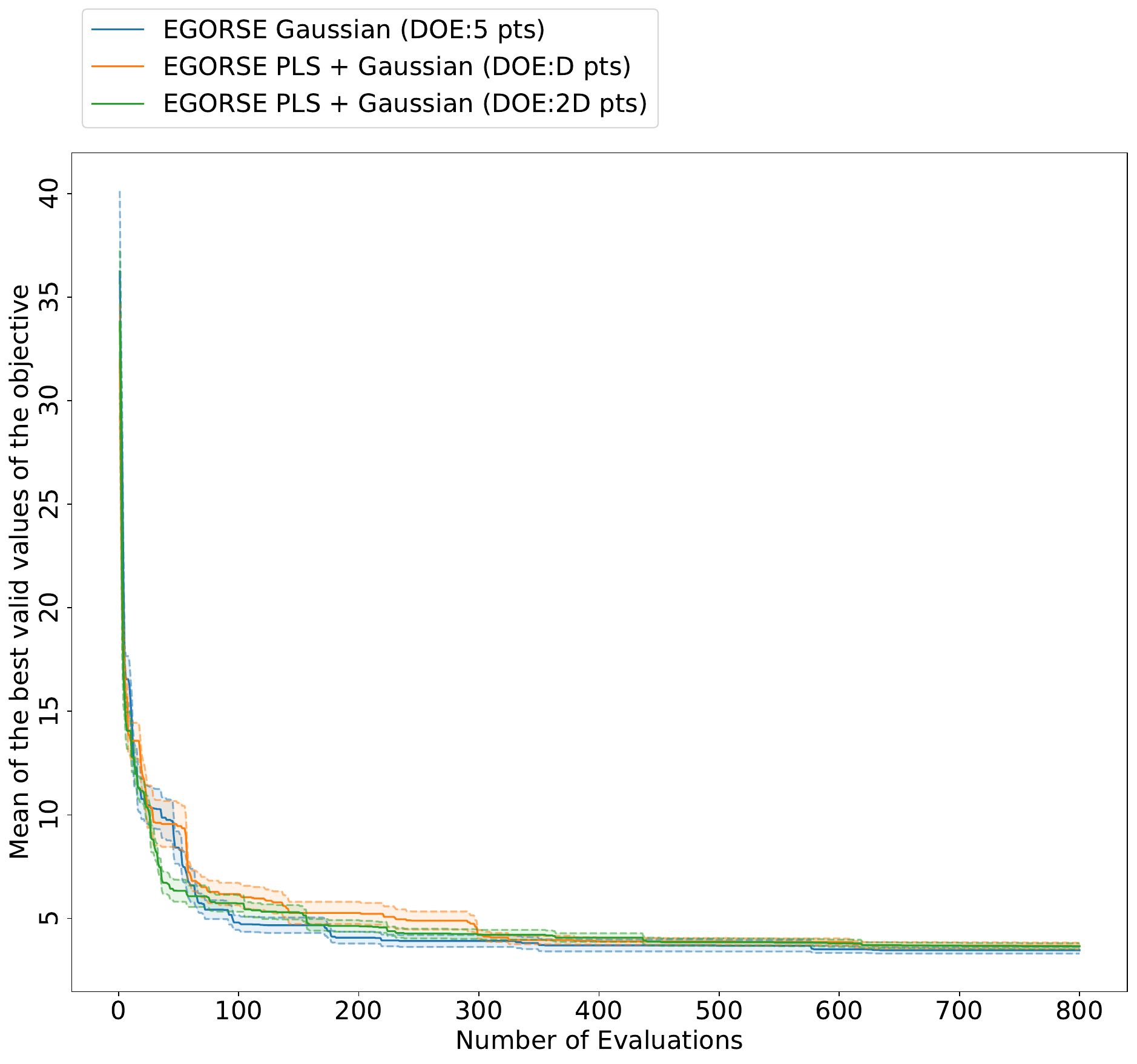}}
    \caption{Convergence plots of 6 versions of \texttt{EGORSE} applied on the MB\_10 problem. The grey vertical line shows the size of the initial DoE.}
    \label{fig:EGORSE_MB10}
\end{figure}

\begin{figure}[!htbp]
    \centering
    \subfloat[DoE of 5 points.\label{fig:EGORSE_MB100_5}]{\includegraphics[width=0.50\textwidth]{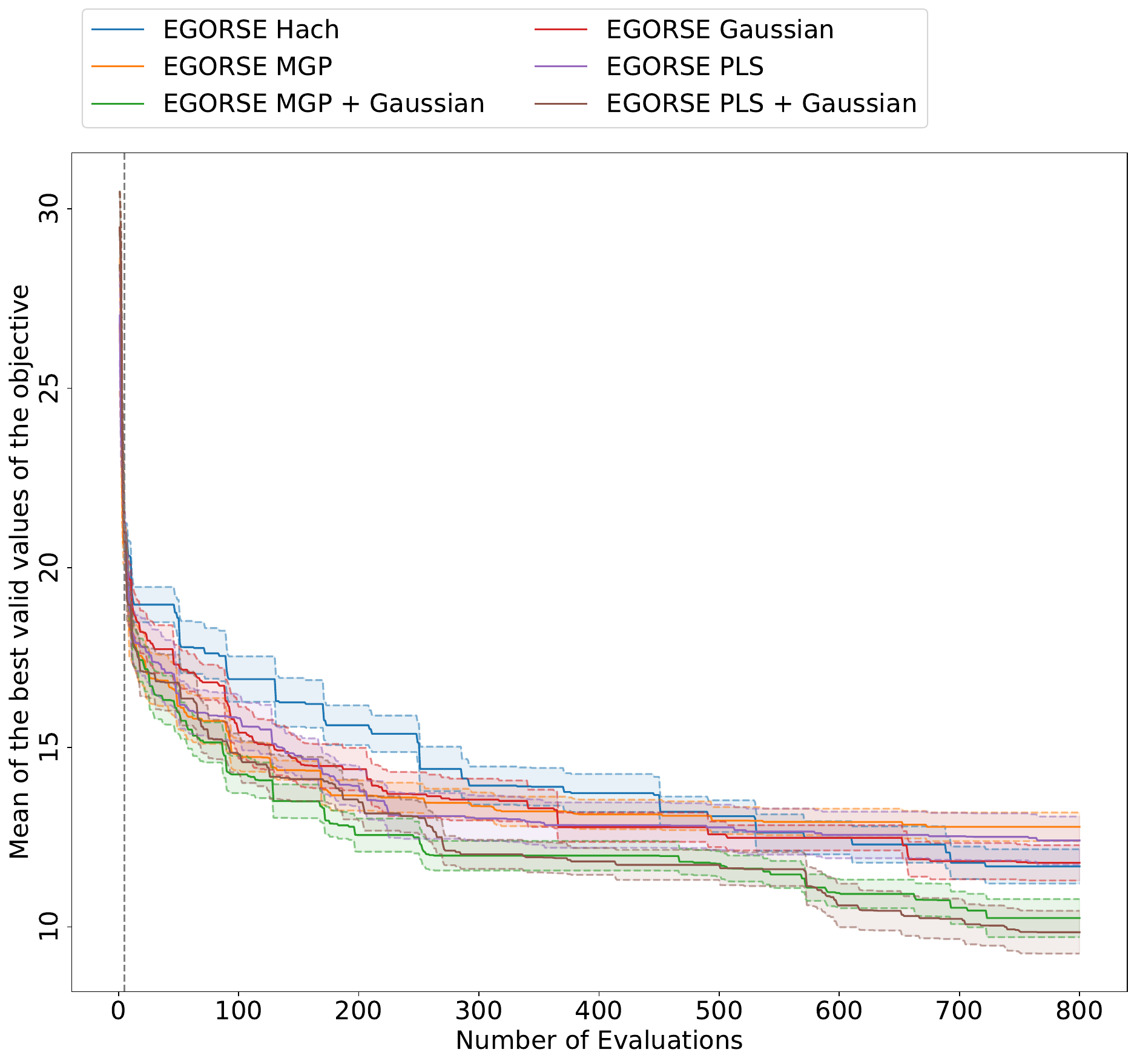}}
    \subfloat[DoE of $d$ points.\label{fig:EGORSE_MB100_d}]{\includegraphics[width=0.50\textwidth]{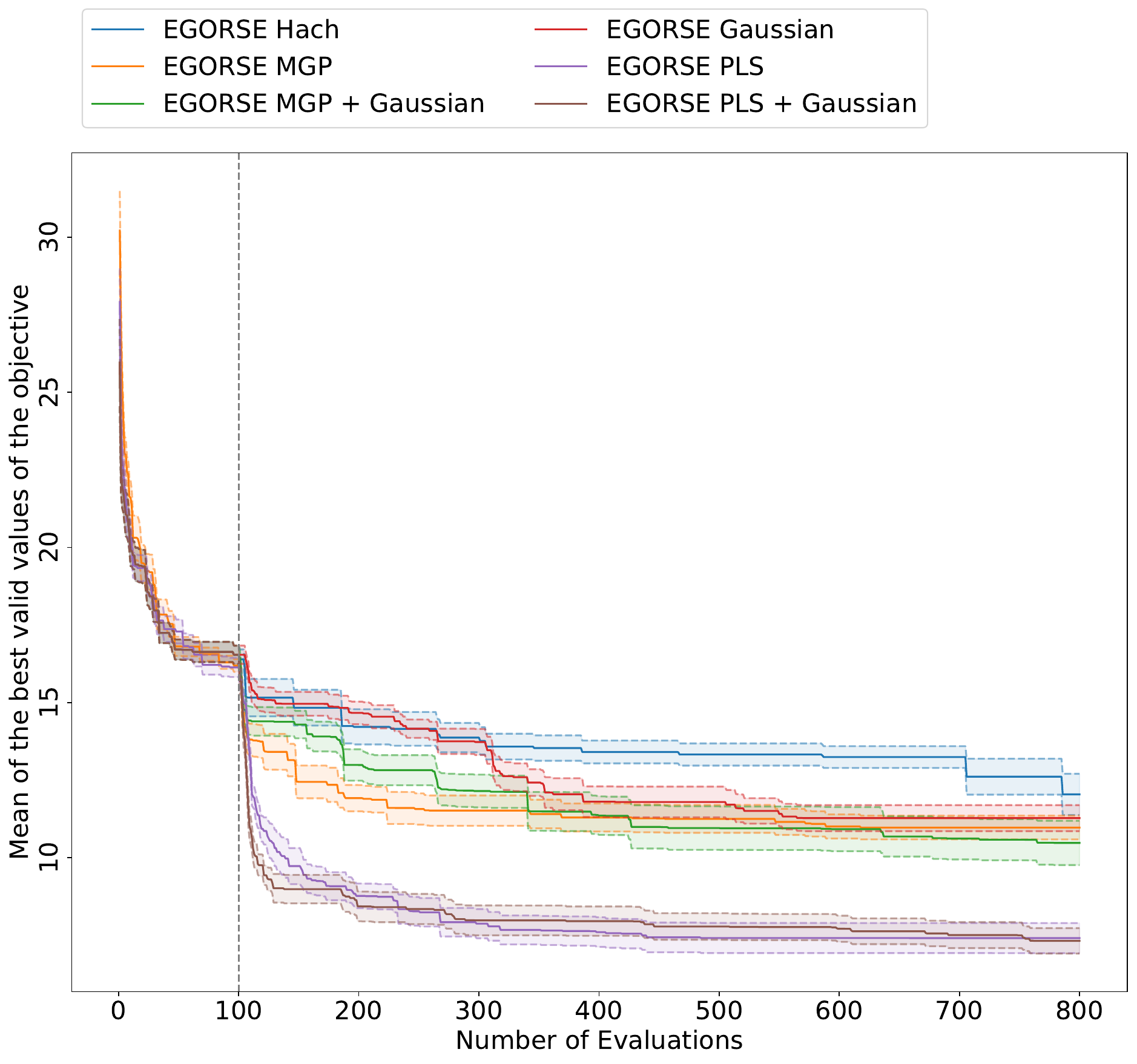}} \\
    \subfloat[DoE of $2d$ points.\label{fig:EGORSE_MB100_2d}]{\includegraphics[width=0.50\textwidth]{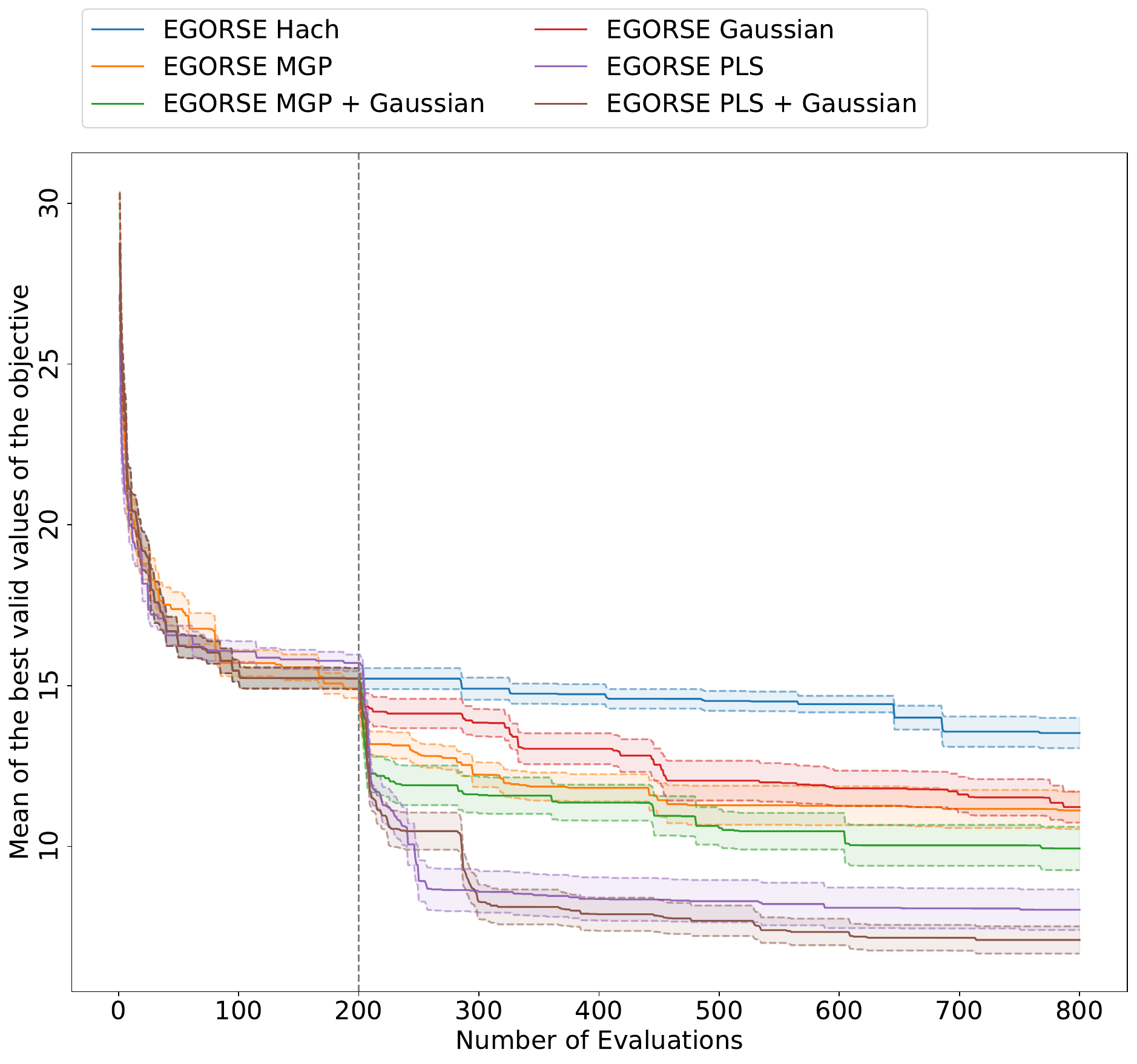}}
    \subfloat[Best versions for each DoE size. \label{fig:EGORSE_MB100_best}]{\includegraphics[width=0.50\textwidth]{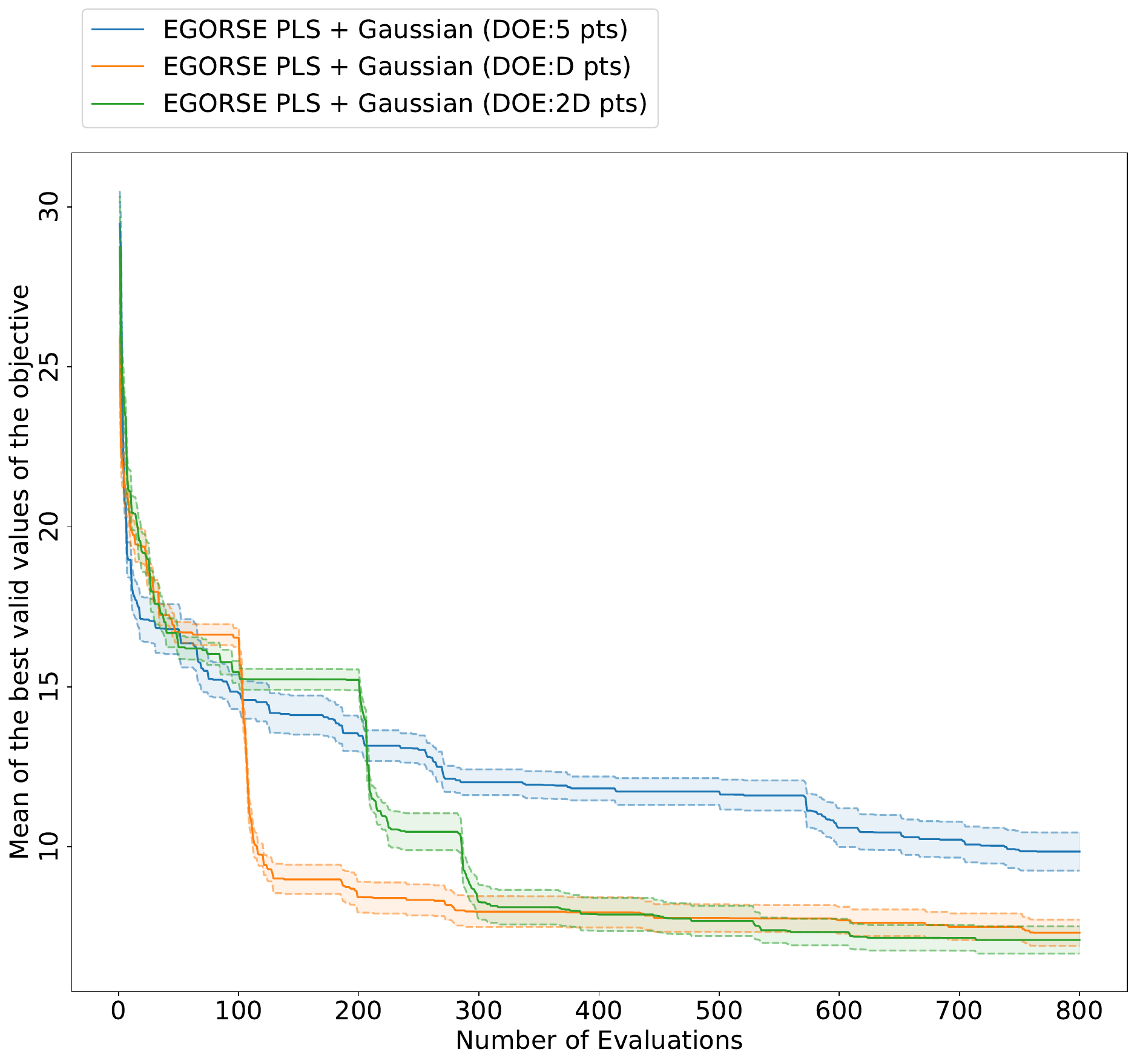}}
    \caption{Convergence plots of 6 versions of \texttt{EGORSE} applied on the MB\_100 problem. The grey vertical line shows the size of the initial DoE.}
    \label{fig:EGORSE_MB100}
\end{figure}

In this section, the convergence robustness and speed of the six \texttt{EGORSE} versions are analyzed with convergence plots (see Section~\ref{sssec:conv}). {\textcolor{black}{Reaching the global minimum (value of 1) of the modified Branin optimization problem in high dimension~\cite{ParrInfillsamplingcriteria2012}  is very hard specially with a limited budget of function evaluations, so the the main purpose here is to compare the convergence speed and the quality of the best value found by  each method.}} 
The convergence plots of the considered \texttt{EGORSE} versions for the MB\_10 problems are displayed in Figure~\ref{fig:EGORSE_MB10}.
In Figure~\ref{fig:EGORSE_MB10_5}, one can see that the \texttt{EGORSE Gaussian} algorithm offers the best performance in term of convergence speed and robustness for an initial DoE of 5 points. 
However, Figures~\ref{fig:EGORSE_MB10_d} and~\ref{fig:EGORSE_MB10_2d} show that \texttt{EGORSE PLS+Gaussian} outperforms \texttt{EGORSE Gaussian} in term of convergence speed and robustness for larger initial DoE. 
Figure~\ref{fig:EGORSE_MB10_best} finally compares these different versions and they are hardly distinguishable for the MB\_10 problem. 
Thus, all the six \texttt{EGORSE} algorithms seem to provide equivalent performance for low dimensional problems.
Figure~\ref{fig:EGORSE_MB100} displays the convergence plots of the six \texttt{EGORSE} versions applied to the MB\_100 problem.
One can see that the different \texttt{EGORSE} versions are much more distinguishable on Figure~\ref{fig:EGORSE_MB100}. 
In fact, Figures~\ref{fig:EGORSE_MB100_5},~\ref{fig:EGORSE_MB100_d} and~\ref{fig:EGORSE_MB100_2d} show that \texttt{EGORSE PLS+Gaussian} provides the best convergence speed and robustness trade-off for all the initial DoE sizes tested.
To select the best number of initial DoE points, Figure~\ref{fig:EGORSE_MB100_best} displays the convergence plots of  \texttt{EGORSE PLS+Gaussian} for the three tested number of points in the initial DoE.
One can easily see that  \texttt{EGORSE PLS+Gaussian} with an initial DoE of $d$ points  provides the best performance in term of convergence speed and robustness.
in terms of evaluations number, the use of an initial DoE of $d$ points allows the algorithm to explore the domain in an interesting direction more quickly than a $2d$ points initial DoE.  
On the contrary, using an initial DoE of 5 points forces the algorithm to seek for the best direction for a long time.

To conclude, choosing the PLS and Gaussian dimension reduction methods with an initial DoE of $d$ points seems the most suitable to obtain the best performance of \texttt{EGORSE}.
Nevertheless, \texttt{EGORSE} capabilities have to be compared with HDBO algorithms to validate its usefulness as it is proposed in the following.

\section{Comparison with state-of-the-art HDBO methods}
\label{sec:num}
\subsection{BO algorithms and setup details}
\label{ssec:setup_all}

\texttt{EGORSE} is now compared to the following state-of-the-art algorithms:
\begin{itemize}
    \item \texttt{TuRBO}~\cite{erikssonScalableGlobalOptimization2019a}: a HDBO algorithm using trust regions to favor the exploitation of the DoE data.
    Tests are performed with the TuRBO\footnote{https://github.com/uber-research/TuRBO}~\cite{erikssonScalableGlobalOptimization2019a} Python toolbox.
    \item   \texttt{EGO-KPLS}~\cite{BouhlelEfficientglobaloptimization2018}: an HDBO method relying on the reduction of the number of GP hyper-parameters.
    This allows to speed up the GP building. 
    The SEGOMOE~\cite{bartoliAdaptiveModelingStrategy2019} Python toolbox is used.
    All the hyper-parameters of this algorithm are the default ones.
    The number of principal components for the KPLS model is set to two.
    \item  \texttt{RREMBO}~\cite{binoisChoiceLowdimensionalDomain2020}: a HDBO method using the random Gaussian transfer matrix to reduce the number of dimensions of the optimization problem through random embedding.
     \texttt{RREMBO}\footnote{https://github.com/mbinois/RRembo} implementation of this algorithm is used.
    The parameters are also set by default.
    \item   \texttt{HESBO}~\cite{nayebiFrameworkBayesianOptimization2019}: a HDBO algorithm using Hash tables to generate the transfer matrix and construct so called hashing-enhanced subspaces. 
    We use the \texttt{HESBO}\footnote{https://github.com/aminnayebi/HesBO} Python toolbox with the default parameters.
\end{itemize}
For \texttt{EGORSE}, the version showing the best performance in term of convergence speed and robustness in Section~\ref{sssec:impact} is selected, i.e.  \texttt{{\texttt{EGORSE}} PLS + Gaussian} with an initial DoE of $d$ points.
To achieve this comparison, $10$ optimizations for each problem and  for each studied method are completed to analyze the statistical behavior of these BO algorithms.
Because of the different strategies implemented in the previously introduced algorithms, a specific test plan must be adopted for each of them.
\begin{itemize}
    \item \texttt{EGORSE}: The test plan of Section~\ref{sssec:conv} is implemented.
    \item  \texttt{TuRBO}: Five trust regions are used with a maximum of $800$ evaluations of the objective function. Note that it is not possible to provide the same initial DoE used for \texttt{EGORSE} in  \texttt{TuRBO}. The number of points generated at the beginning of the algorithm is thus imposed to $d$. The EI acquisition function is chosen.
     \item  \texttt{EGO-KPLS}: The optimization is performed with a maximum evaluation number of $800$, with the initial DoE used in \texttt{EGORSE} and with the EI acquisition function. 
     \item  \texttt{RREMBO \& HESBO}: $20$ optimizations of $20d_e$ evaluations are performed for each of the initial DoE for  \texttt{RREMBO} and  \texttt{HESBO}. These $20$ optimizations are then concatenated and considered as a unique optimization process. The number of effective directions are imposed to $d_e=2$ and the acquisition function is EI.  \texttt{RREMBO} and  \texttt{HESBO} are equivalent to \texttt{EGORSE} without supervised dimension reduction method.
\end{itemize}
\subsection{Results analysis}
\label{ssec:bo_mb:res}

\begin{figure}[!hbtp]
    \centering
    \subfloat[MB\_10 iteration convergence plot. \label{fig:ALL:mb10_nb}]{\includegraphics[width=0.50\textwidth]{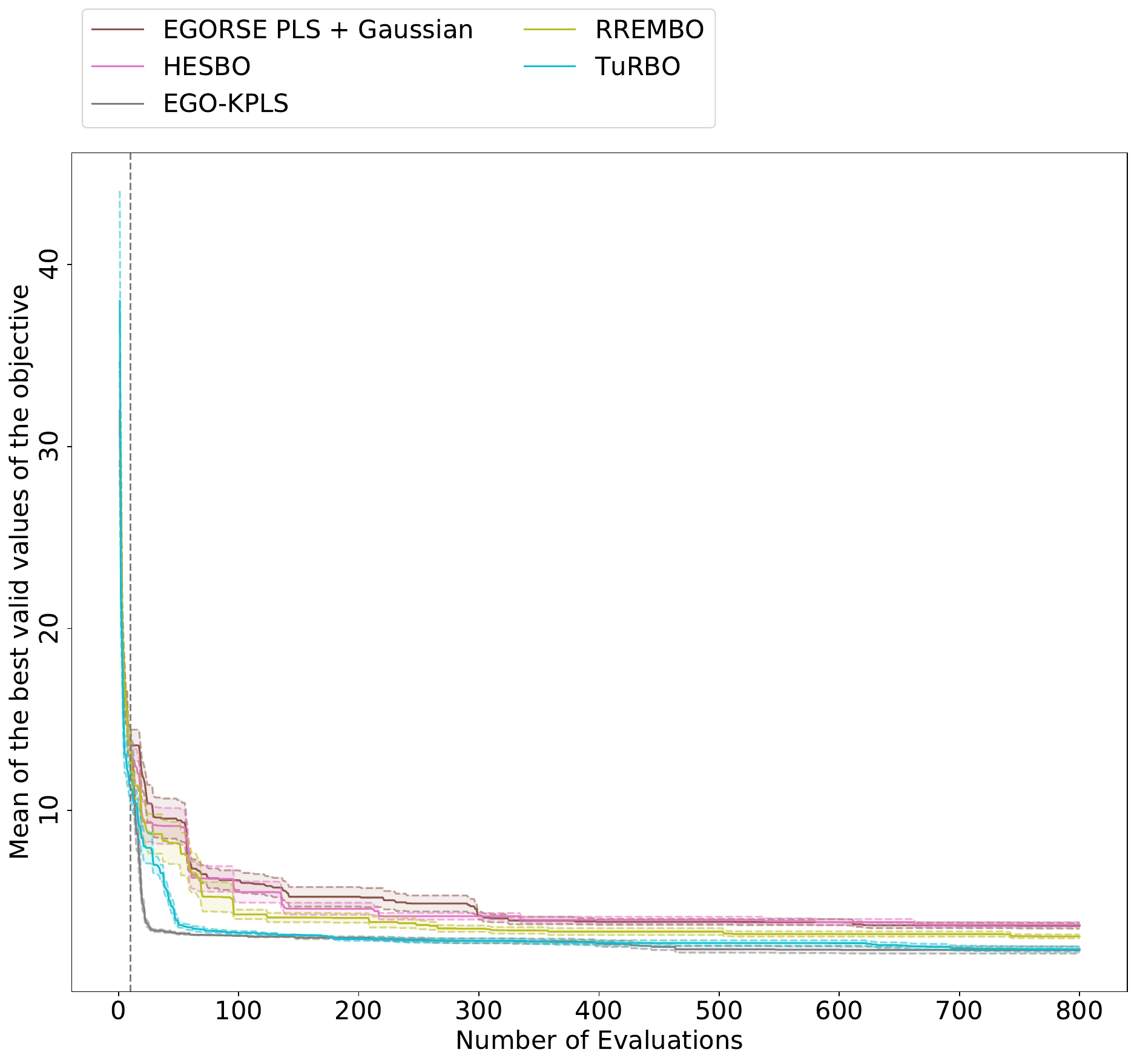}}
    \subfloat[MB\_100 iteration convergence plot. \label{fig:ALL:mb100_nb}]{\includegraphics[width=0.50\textwidth]{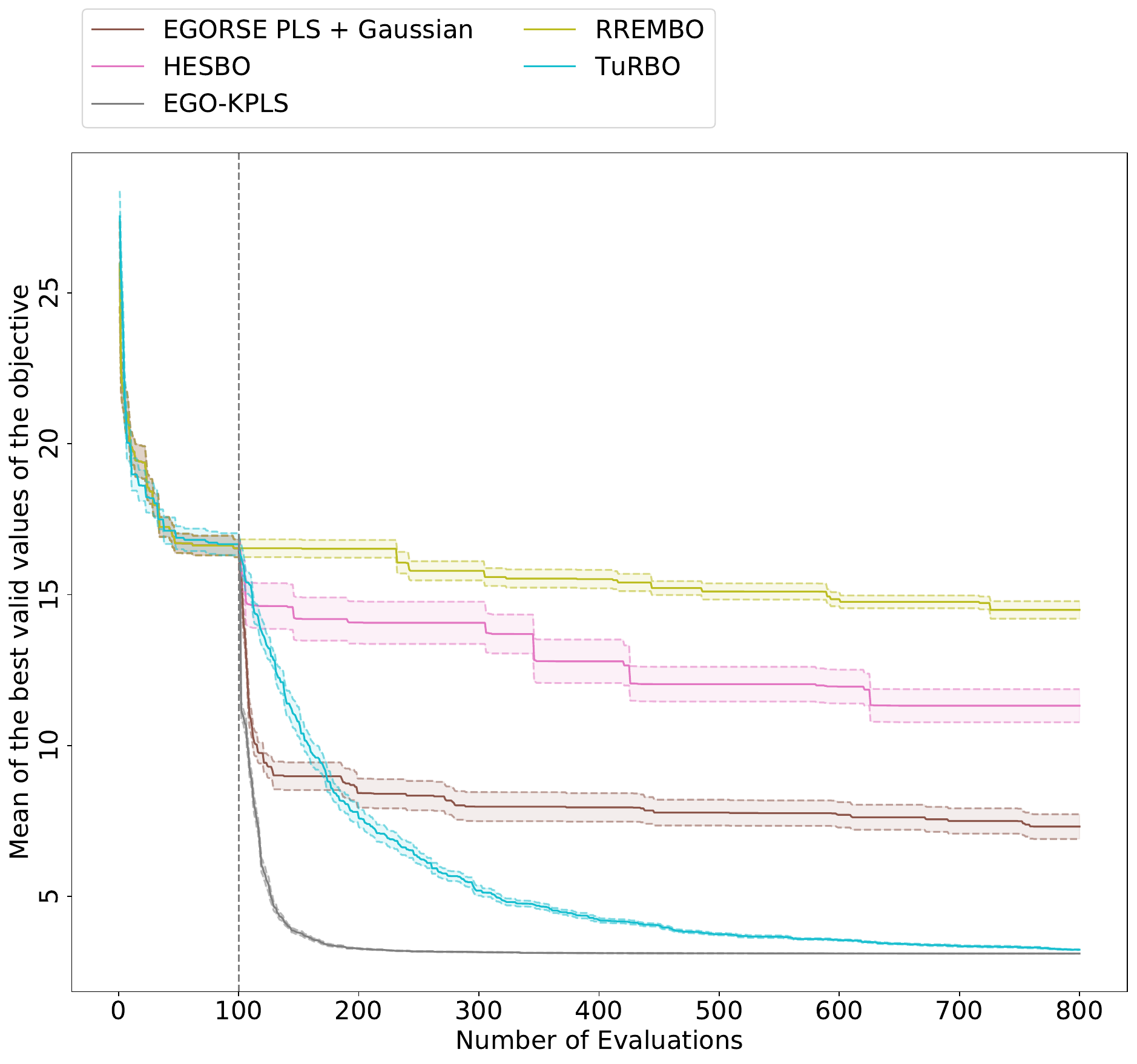}} \\
    \subfloat[MB\_10 time convergence plot. \label{fig:ALL:mb10_time}]{\includegraphics[width=0.50\textwidth]{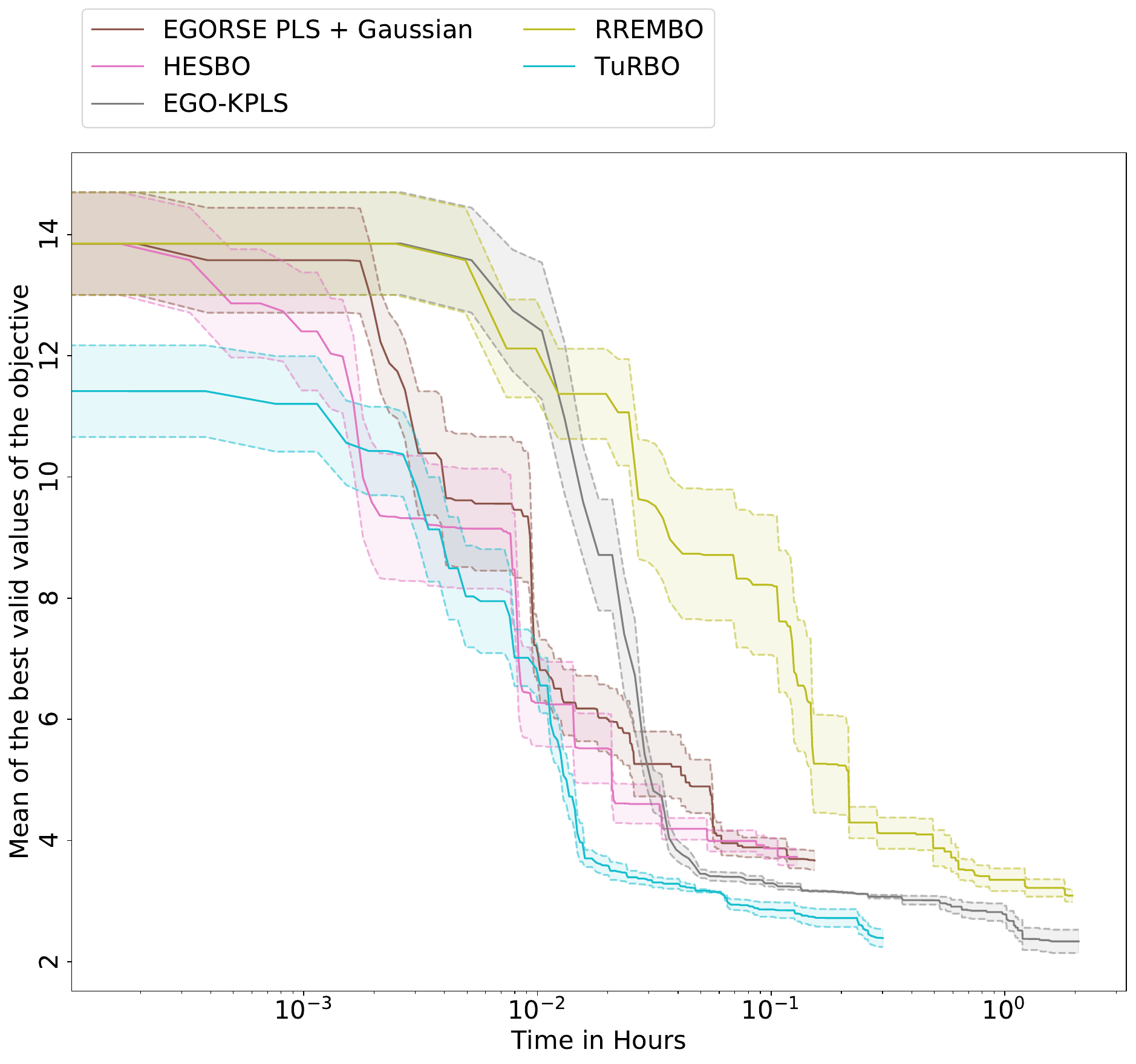}}
    \subfloat[MB\_100 time convergence plot. \label{fig:ALL:mb100_time}]{\includegraphics[width=0.50\textwidth]{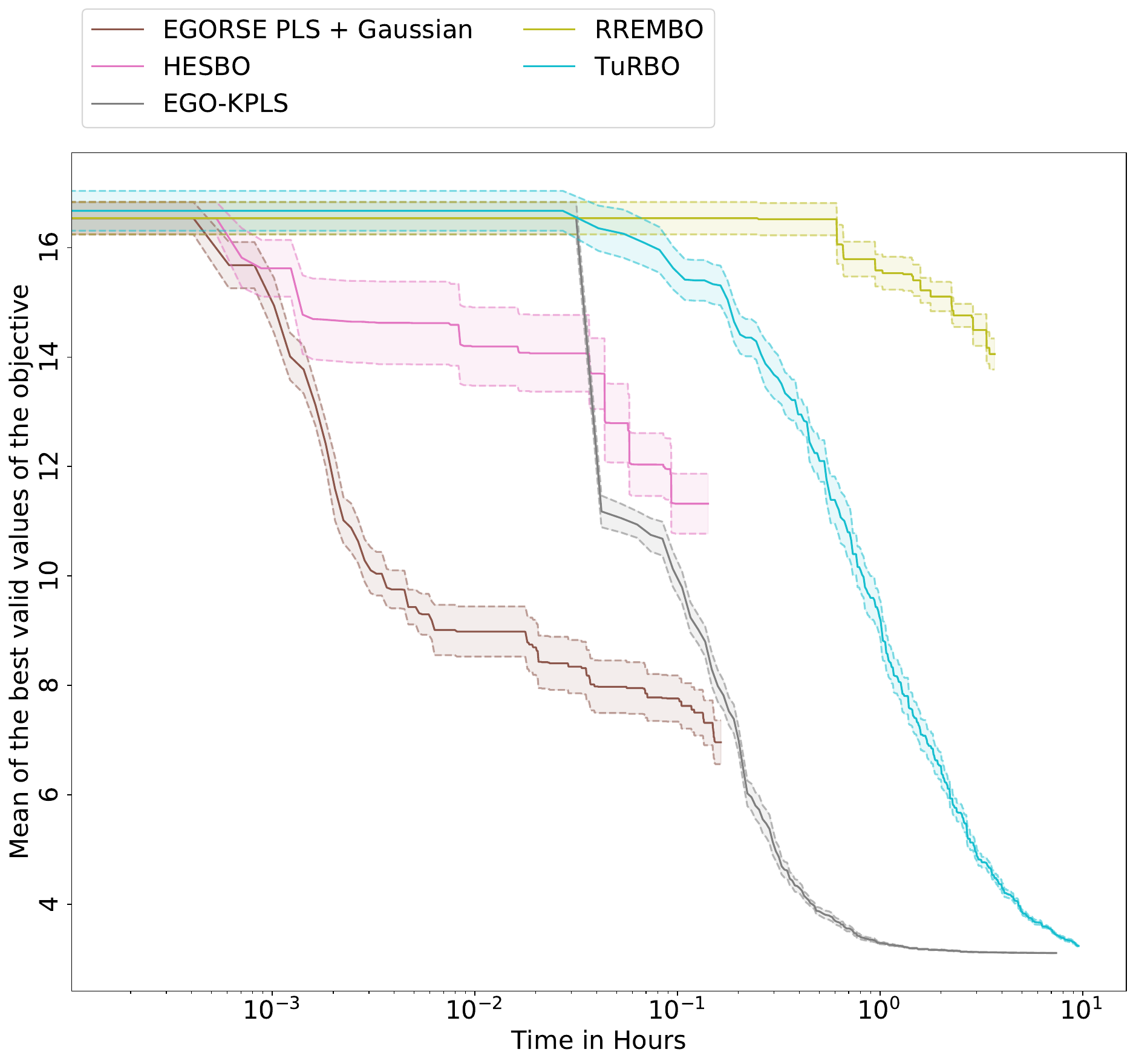}}
    \caption{Iteration and time convergence plots for 5 HDBO algorithms on the MB\_10 and MB\_100 problems. The grey vertical line shows the size of the initial DoE.}
    \label{fig:ALL}
\end{figure}

In this section, a comparison between \texttt{EGORSE} and the four studied algorithms is performed in term of robustness, convergence speed both in CPU time and in number of iterations.
Figure~\ref{fig:ALL} provides the iteration convergence plots, as introduced in Section~\ref{sssec:conv}, and the time convergence plots drawing the evolution of the means of the best discovered function values against the CPU time. 

Figure~\ref{fig:ALL:mb10_nb} shows that  \texttt{TuRBO} and  \texttt{EGO-KPLS} are converging the fastest and with a low standard deviation.
Moreover, the convergence plots of  \texttt{EGORSE},  \texttt{RREMBO} and  \texttt{HESBO} are hardly distinguishable. 
Figure~\ref{fig:ALL:mb100_nb} displays that  \texttt{EGO-KPLS} converges the fastest to the lowest values with a low standard deviation.
 \texttt{TuRBO} is also providing good performance even if it converges slower than  \texttt{EGO-KPLS}. 
Regarding the three methods using dimension reduction procedure,  \texttt{EGORSE} is converging to the lowest value with a relatively low standard deviation.
The good performance of  \texttt{EGO-KPLS} and  \texttt{TuRBO} is certainly due to the ability of these algorithms to search all over $\Omega$, which is not the case for other methods.
However, when the dimension of $\Omega$ increases, the ability to search all over $\Omega$ becomes a drawback. In fact, a complete search in $\Omega$ is intractable in time is this case.

Figures~\ref{fig:ALL:mb10_time} and~\ref{fig:ALL:mb100_time} depict the convergence CPU time necessary to obtain the regarded value.
First, the RREMBO, TuRBO and EGO-KPLS complete the optimization procedure in more than 8 hours on the MB\_100 problem against an hour on the MB\_10 problem.
This suggests that  \texttt{RREMBO},  \texttt{TuRBO} and  \texttt{EGO-KPLS} are intractable in time for larger problems.
Then, one can easily see that \texttt{EGORSE} is converging the fastest in CPU time than the other algorithms on the MB\_100 problem.
In fact, the computation time needed to find the enrichment point is much lower than the one for  \texttt{TuRBO},  \texttt{EGO-KPLS} and  \texttt{RREMBO}.
This was sought in the definition of the enrichment sub-problem introduced in Section~\ref{ssec:sub_prob}.
Finally,  \texttt{EGORSE} is converging to a lower value than  \texttt{HESBO} in a similar amount of time. 
Thus,  \texttt{EGORSE} seems more interesting to solve HDBO problems than the studied algorithms.
Note that only  \texttt{HESBO} and  \texttt{EGORSE} are able to perform an optimization procedure on HDBO problems.

To conclude this section, several sets of hyper-parameters of  \texttt{EGORSE} have been tested on two problems of dimension 10 and 100. 
The analysis of the obtained results has shown that  \texttt{EGORSE PLS+Gaussian} with an initial DoE of $d$ points is performing the best. 
A comparison of  \texttt{EGORSE} with HDBO algorithms has also been carried out. 
It has pointed out that  \texttt{TuRBO},  \texttt{EGO-KPLS} and  \texttt{RREMBO} are intractable in time for HDBO problems.
Furthermore,  \texttt{EGORSE} has appeared to be the most suitable to solve HDBO problem efficiently. 

\section{Evaluation of \texttt{EGORSE} on a high dimensional planning optimization}
\label{sec:rover}

In this section, the \texttt{EGORSE} algorithm is evaluated to find an optimal path planning problem using 600 design optimization variables. 
\subsection{Problem definition and implementation details}

\begin{figure}[!htp]
    \centering
    \includegraphics[width=0.8\textwidth]{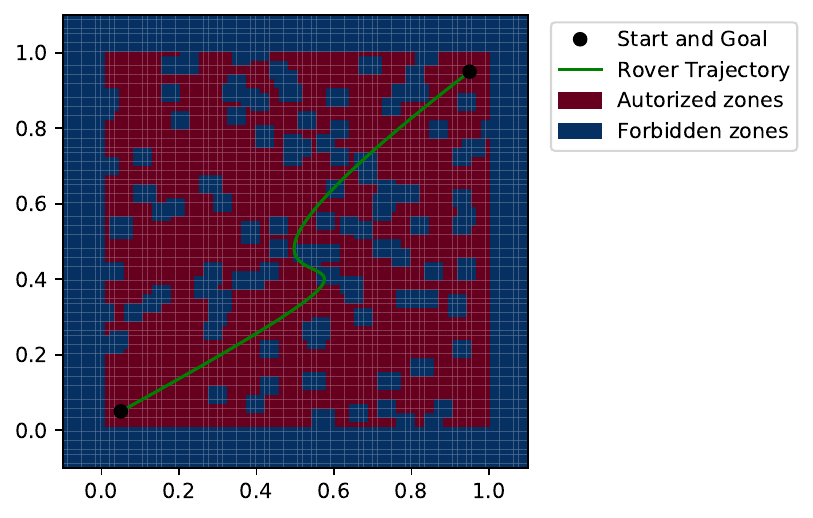}
    \caption{Example of a robot trajectory in a forest.}
    \label{fig:rover_dom}
\end{figure}

The Rover\_600 path planning problems relies on the same idea than MB\_d problems except that the objective function is a adjustment of the Rover\_60~\cite{WangBatchedHighdimensionalBayesian2018} problem.
It consists on a robot routing from a starting point $x_{start}$ to a goal point $x_{goal}$ in a forest. 
The robot trajectory is a spline defined by 30 control points. 
These points, including the starting and the goal ones, are the design variables of the optimization problem that belong to $\Omega=[0,1]^{60}$.
The objective function is minimal when the robot follows the shortest trajectory without meeting a tree.
The minimum of the function is $f_{min} = -5$.
Figure~\ref{fig:rover_dom} gives an example of trajectory.

To increase the number of design variables, the problem is normalized in $\Omega=[-1,1]^{60}$, a random matrix $\A_d \in \mathbb{R}^{60  \times 600}$ is generated such that all $\x \in [-1,1]^d$, $\A_d \x = \bm{u} \in [-1,1]^{60}$.
An objective function Rover\_600, where $d=600$ is the number of design variables, is defined such that $\text{Rover\_600}(\x) = \text{Rover\_60}(\A_d\x)$.
Eventually, we solve the following optimization problem:
\begin{equation}
    \min\limits_{x \in [-1,1]^600} \text{Rover\_600}(\x) = \min\limits_{x \in [-1,1]^{60}} \text{Rover\_60}(\A_d\x).
\end{equation}

We note that, in our tests, \texttt{TuRBO},  \texttt{EGO-KPLS} and  \texttt{RREMBO} algorithms are not included anymore in this comparison as they are intractable in time for optimization problems with more than a hundred design variables (see Section~\ref{ssec:bo_mb:res}). So, only \texttt{EGORSE PLS+Gaussian} is compared to \texttt{HESBO}. The same test plan as in Section~\ref{ssec:setup_all} is used with $200$ optimizations (instead of $20$) of $20d_e$ evaluations.

\subsection{Result analysis}
\begin{figure}[!hbt]
    \centering
    \subfloat[Time convergence plot. \label{fig:Rover:Time}]{\includegraphics[width=0.50\textwidth]{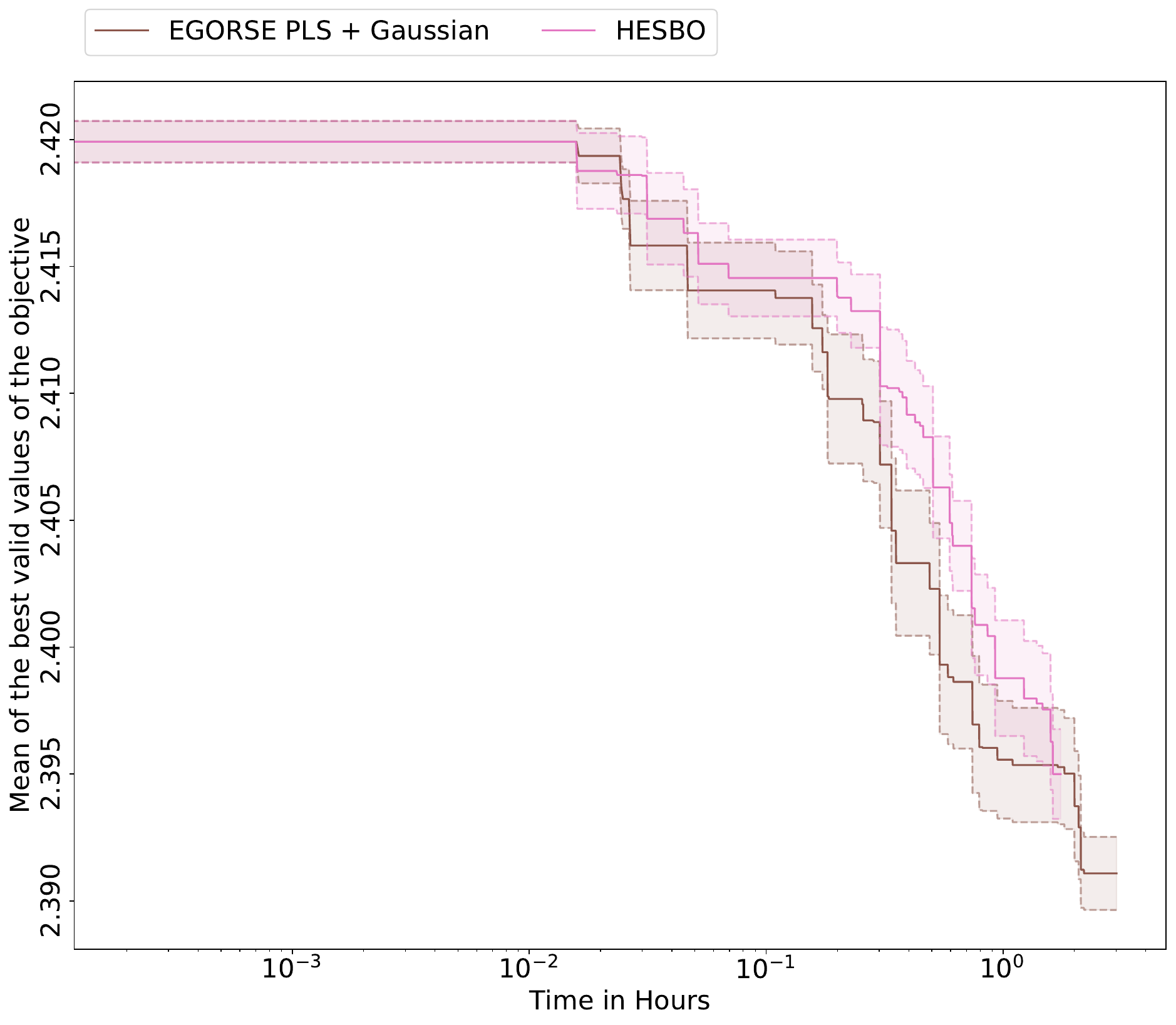}}
    \subfloat[Iteration convergence plot. \label{fig:Rover:nbCall}]{\includegraphics[width=0.50\textwidth]{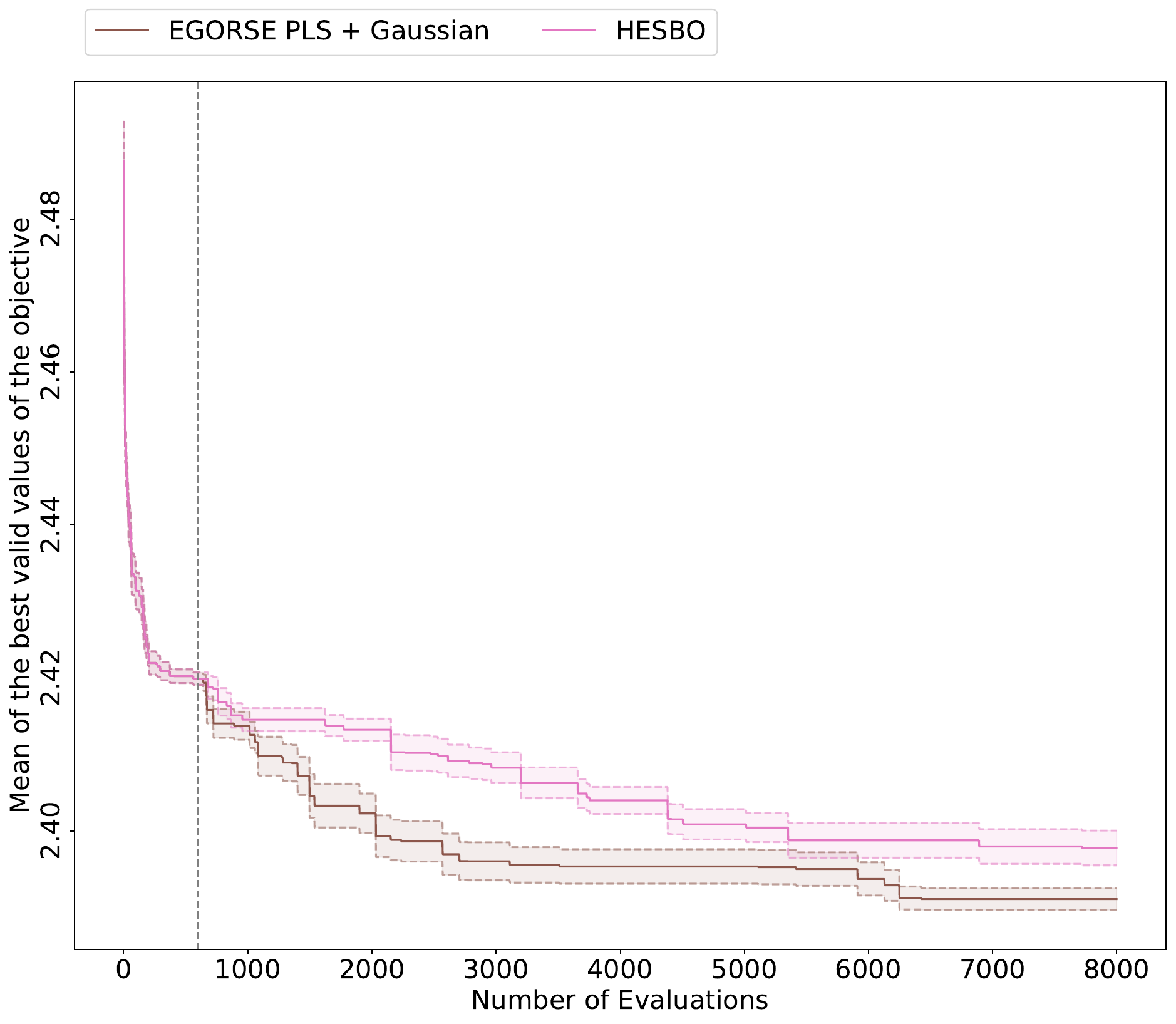}}
    \caption{CPU time convergence plots of  \texttt{HESBO, RREMBO, TuRBO, EGO-KPSL} on the Rover\_600 problem with an initial DoE of $d$ points. The grey vertical line shows the size of the initial DoE.}
    \label{fig:Rover}
\end{figure}

Time and iteration convergence plots of  \texttt{EGORSE} and  \texttt{HESBO} are depicted in Figure~\ref{fig:Rover}.
Figure~\ref{fig:Rover:nbCall} clearly shows that  \texttt{EGORSE} is converging fast to the lowest objective value with  a very low standard deviation. 
However, the obtained objective value is larger than the known optimal one (i.e. $f_{min}=-5$).
This is due to two main reasons that may be addressed in further research:
\begin{itemize}
    \item The number of effective directions used in  \texttt{EGORSE} (i.e. $d_e=2$) is much lower than the actual number of effective directions (i.e. $d_e=60$). The effective search space is not covering the space in which Rover\_600 is varying.
    \item The dimension reduction method PLS is global. The local variations of the function, in which the global optimum can be located, are thus deleted. Even if this problem is tackled by searching in randomly generated subspace,  \texttt{EGORSE} cannot provide better results. 
\end{itemize}
Figure~\ref{fig:Rover:Time} shows that  \texttt{HESBO} is performing the optimization procedure faster than  \texttt{EGORSE}.
In fact,  \texttt{HESBO} does not solve any quadratic problem at each iteration.
However, one can see that the time difference is not significant.

\section{Conclusion}
\label{sec:conclusion}

This paper introduces EGORSE, a high-dimensional efficient global optimization using both random and supervised embeddings to tackle expensive to compute black-box optimization problems.
EGORSE shows a high potential to tackle the two main drawbacks of most existing HDBO methods with a very competitive computational effort and with an accurate definition of the new bounds related to the reduced optimization problem.
A parametric study on the hyper-parameters of EGORSE has shown that combining both PLS and the random Gaussian reduction methods provides the best results. On very large-scale optimization problem, i.e. d=600, EGORSE gives a lower minimal value than HESBO in an equivalent computational time. 
Thus EGORSE outperforms by far all the state-of-the-art HDBO solvers.
Nonetheless, EGORSE is still unable to find the global optimum in many cases as the effective dimension of the problem is much larger than the effective dimension used in the algorithm. {\color{black}Overall, EGORSE shows very good performance compared to state-of-the-art high-dimensional BO solvers.}


{\color{black} We believe that the results from EGORSE encourage further study and enhancements to tackle realistic aerospace engineering optimization problems. In fact, all realistic aerospace optimization problems involve constraints~\cite{priemEfficientApplicationBayesian2020a, raponi2020high,prestch2024}. Extending EGORSE to solve high-dimensional constrained black-box optimization problems is necessary to address issues where constraints are typically related to the amount of oil, the range of an aircraft, or the size of its wings. Another important perspective is related to the asymptotic convergence properties of EGORSE. Convergence guarantees can be obtained, for instance, by incorporating trust-region techniques~\cite{YDIOUANE_2023}. This could help EGORSE focus on promising regions, improve its convergence rate to local optima, and use a stopping criterion based on a stationarity measure rather than a fixed budget. 

One current limitation of EGORSE is related to the effective dimension (i.e., the dimension of the reduced space), which needs to be specified by the user. However, in real-world applications, users do not necessarily know the effective dimensions. In this context, an automatic estimation of the number of effective dimensions, inspired by \cite{SciTech_cat}, would be of high interest toward developing an effective dimension-free methodology.}




\section*{Acknowledgments}
This work is part of the activities of ONERA - ISAE - ENAC joint research group.
The authors would like to express their special thanks to Nina Moello for her test on the implementation of the log-likelihood gradient and Hessian in the opensource Python toolbox SMT\footnote{https://github.com/SMTorg/smt} and for her implementation of the HESBO algorithm in our computer environment.

\bibliography{main}

\end{document}